\documentclass[a4paper,a4wide,11pt]{article}
\usepackage{a4wide}
\usepackage{theorem}
\usepackage{amsmath}
\usepackage{array}
\usepackage{amssymb}
\usepackage{amsfonts}
\usepackage[french]{babel}
\usepackage{epsf}
\usepackage{epsfig}

\newtheorem{theo}{\indent Theorem\newline}[section]

{\theorembodyfont{\rmfamily}%
\newtheorem{rem}[theo]{\noindent Remark}}
{\theorembodyfont{\rmfamily}%
 \theoremstyle{break}%
\newtheorem{reml}[theo]{\noindent Remark}}
\newtheorem{prop}[theo]{\indent Proposition\newline}
\newtheorem{lemma}[theo]{\indent Lemma\newline}
\newtheorem{cor}[theo]{\indent Corollary\newline}

 \def\N{{\Bbb{N}}}
\def\Z{{\Bbb{Z}}}

\def\R{{\Bbb{R}}}
\def\C{{\Bbb{C}}}

\newcommand{\Div}{\mathop{\rm div}\nolimits}
\newcommand{\Jac}{\mathop{\rm Jac}\nolimits}
\newcommand{\Pic}{\mathop{\rm Pic}\nolimits}

\newlength{\indentation}%
\makeatletter
\newcommand\@makefntextsans[1]{%
    \parindent 0em%
    \noindent%
    \hb@xt@0em{\hss}%
    #1}
\def\footnotetextsans{%
     \@ifnextchar [\@xfootnotenextsans%
       {\@footnotetextsans}}
\def\@xfootnotenextsans[#1]{%
  \begingroup%
     \csname c@\@mpfn\endcsname #1\relax%
  \endgroup%
  \@footnotetextsans}
\long\def\@footnotetextsans#1{\insert\footins{%
    \reset@font\footnotesize%
    \interlinepenalty\interfootnotelinepenalty%
    \splittopskip\footnotesep%
    \splitmaxdepth \dp\strutbox \floatingpenalty \@MM%
    \hsize\columnwidth \@parboxrestore%
    \color@begingroup%
      \@makefntextsans{%
        \rule\z@\footnotesep\ignorespaces#1\@finalstrut\strutbox}
    \color@endgroup}}
\makeatother

\newlength{\avantsectionapp}%
\newlength{\apressectionapp}%
\newlength{\avantsubsectionapp}%
\newlength{\apressubsectionapp}%
\newcounter{sectionapp}%
\newcounter{subsectionapp}[sectionapp]%
\begin{document}

\cleardoublepage
\title{Real structures on minimal ruled surfaces}
\author{Jean-Yves Welschinger \\}
\maketitle

\makeatletter\renewcommand{\@makefnmark}{}\makeatother
\footnotetextsans{Keywords :  Ruled surface, real algebraic surface.}
\footnotetextsans{Classification AMS : 14J26, 14P25.}

{\bf Abstract :}

In this paper, we give a complete description of the deformation
classes of real structures on minimal ruled surfaces. In particular, we show 
that these classes
are determined by the topology of the real structure, which means,
using the terminology of \cite{KhDg2}, that real minimal ruled surfaces are
quasi-simple. 
As an intermediate result, we obtain the
classification, up to conjugation, of real structures on decomposable
ruled surfaces.

\section*{Introduction}

Let $X$ be a smooth compact complex surface. A {\it real structure} 
on $X$
is an antiholomorphic involution $c_X : X \to X$. The {\it real part} of
$(X , c_X)$ is by definition the fixed point set of $c_X$. If $X$ admits
a holomorphic submersion on a smooth compact complex irreducible curve
$B$ whose fibers 
have genus zero, then it is called a {\it minimal ruled surface}. These
surfaces are all algebraic, minimal and of kodaira dimension $- \infty$
 (see \cite{Beau}). Real minimal ruled surfaces are one of the few 
examples of real algebraic surfaces of special type whose classification under
real deformation is not known, see the recent results 
\cite{KhDg2}, \cite{KhDg}, \cite{Cat} and the survey
\cite{Kh} for detailed history and references. The purpose of this paper
is to fill this gap. Since all the ruled surfaces considered in 
this paper will be
minimal, from now on we will call them ``ruled'' rather than ``minimal ruled''.

Rational surfaces are well known
(see \cite{KhDg2}), so we can restrict ourselves to non-rational ruled
surfaces. The ruling $p : X \to B$ is then unique and any real structure
$c_X$ on $X$ is fibered over a real structure $c_B$ on $B$ in the sense that
$c_B \circ p = p \circ c_X$. The topology of the real part of $X$ 
as well as the topology of the real curve $(B , c_B)$ provide us with a 
topological invariant under real deformation which we call
{\it the topological type} of the surface. This invariant is encoded by a quintuple
 of integers : the number of tori and Klein bottles of $\R X$, the genus of $B$, 
the number of components of $\R B$ and the type of $(B , c_B)$ 
(see \S \ref{subsectiontoptype}). The main result of this paper is
the following (see theorem
\ref{theoremdeformation} and proposition \ref{proposexistence}) :

\begin{theo}
\label{theointro}
Two real (minimal) non-rational ruled surfaces are in the same real deformation
class if and only if they have the same topological type and homeomorphic
quotients. Moreover, any allowable quintuple of integers is realized as the topological
type of a real non-rational ruled surface.
\end{theo}

Note that as soon as the bases of the surfaces have non-empty real parts,
the condition on the quotients can be removed. A quintuple of 
integers
is called {\it allowable} when it satisfies the few obvious conditions
satisfied by topological types of real non-rational ruled surfaces,
see \S \ref{subsectiontoptype} for a precise definition.
Remember that any compact complex surface lying in the deformation
class of a non-rational ruled surface is itself a non-rational
 ruled surface (see \cite{BPV} for example). A definition of real deformation
classes can be given as follows. 
Equip the Poincar\'e's disk $\Delta \subset \C$ with the complex
conjugation $conj$. A {\it real deformation} of surfaces is a proper
holomorphic submersion $\pi : Y \to \Delta$ where $Y$ is a complex manifold 
of dimension $3$ equipped with a real structure $c_Y$ and $\pi$ satisfies 
$\pi \circ c_Y =
conj \circ \pi$. Then, when $t \in ]-1 , 1 [ \in \Delta$, the fibers $Y_t =
\pi^{-1} (t)$ are invariant under $c_Y$ and hence are compact real
surfaces. Two real surfaces $X'$ and $X''$ are said
to be 
{\it in the same deformation class} if there exists a chain $X'=X_0,
\dots , X_k=X''$ of compact real surfaces such that for every 
$i \in \{ 0 , \dots , k-1 \}$, the
surfaces $X_i$ and $X_{i+1}$ are isomorphic to some real fibers of a
real deformation.

Remember that every ruled surface is the 
projectivisation $P (E)$ of a rank two complex vector bundle $E$ over $B$ (see
\cite{Beau}). Moreover $P (E)$ and $P (E')$ are isomorphic if and only
if $E' = E \otimes L$ where $L$ is a complex line bundle over $B$. A
ruled surface is said to be {\it decomposable} if $E$ is decomposable,
that is if $E$ is the direct sum of two complex line bundles.
The paper is organized as follows. In the first section, we give 
constructions of some particular real structures on decomposable ruled
surfaces. In the second section we obtain a
classification, up to conjugation, of real structures on decomposable
ruled surfaces (see theorem \ref{theoremrealstruct}). This result,
of independant interest, plays a crucial r\^ole in the proof of theorem
\ref{theointro}. In this section is also given a result independant of real 
algebraic geometry, which concerns the
lifting of automorphisms of the ruled surface $X$ to automorphisms of
the rank two vector bundle $E$, see proposition \ref{propaut}. Finally, the
third section is devoted to the proof of theorem \ref{theointro}. This gives
a complete description of the deformation
classes of real structures on ruled surfaces. In particular, it shows
 that these classes
are determined by the topology of the real structure, which means,
using the terminology of \cite{KhDg2}, that real ruled surfaces are
quasi-simple.

\section{Construction of some particular real structures}

\subsection{Meromorphic functions and real structures}

Let $B$ be a smooth compact complex irreducible curve. Denote by $\Pic(B)$ 
the group of 
complex line bundles over $B$. This group is identified with the group of
divisors modulo principal ones. Let $\phi : B \to B$ be a holomorphic
or anti-holomorphic automorphism, and let
$D=\sum_{i=1}^k n_i p_i$, $n_i \in \Z$, $p_i \in B$, be a divisor on
$B$. Then we denote by $\phi^* (D)$ the divisor $\sum_{i=1}^k n_i \phi^{-1} 
(p_i)$ and by $\phi (D)$ the divisor $\sum_{i=1}^k n_i \phi 
(p_i)$. The morphism on the 
quotient $ \Pic(B)$ of the group of divisors induced by $\phi^*$ will
also be denoted by $\phi^*$. We denote by $L_0$ the trivial line bundle
over $B$ and by $L^*$ the line bundle dual to $L$, so that $L \otimes
L^* = L_0$.

Suppose from now on that $B$ is equipped with a {\it real structure} $c_B$,
that is an anti-holomorphic involution $c_B$.

\begin{lemma}
\label{lemmaL}
Let $L \in \Pic (B)$ be a line bundle such that $c_B^* (L)=L$. Then,
for every divisor $D$ associated to $L$, there exists a meromorphic
function $f_D$ on $B$ such that $\Div (f_D) = c_B (D) - D$ and $f_D
\times \overline{f_D \circ c_B} = \pm 1$. 
\end{lemma}

{\bf Proof :}

By assumption, $D$ and $c_B (D)$ are linearly equivalent. As a
consequence, there exists a meromorphic function $f$ such that $\Div
(f) = c_B (D) - D$. Then, $g = \overline{f \circ c_B}$ is a
meromorphic function on $B$ satisfying $\Div (g) = D -c_B (D) $. So
$fg$ is a holomorphic function on $B$. This means 
that there exists
a constant $\lambda \in \C^*$ such that  $f
\times \overline{f \circ c_B} = \lambda$.

But for all $x \in B$, 
$$\lambda = f
\times \overline{f \circ c_B} (c_B (x))= f \circ c_B (x) \times \overline{f
(x)} = \overline{ f(x) \times \overline{f \circ c_B (x)}} =
\overline{\lambda }$$ 

Thus $\lambda \in \R^*$, and we define $f_D =
\frac{1}{\sqrt{|\lambda|}} f$. $\square$\\

\begin{reml}
\label{remarkL}

As soon as $\R B$ is non-empty, $f_D
\times \overline{f_D \circ c_B} = +1$, since for every $x \in \R B$ we
have $f_D
\times \overline{f_D \circ c_B} (x) = |f (x)|^2 \geq 0$. Nevertheless,
when $\R B = \emptyset$, there always exists a divisor $D$ on $B$, of
degree congruent to $g(B) -1 \mod (2)$ where $g(B)$ is the genus of
$B$, such that $f_D
\times \overline{f_D \circ c_B} = -1$ (see \cite{GH}, proposition
2.2). 
Note also that the  function $f_D$ given by lemma \ref{lemmaL}
is not unique, since for 
every constant $\lambda \in \C$ such that $|\lambda| = 1$, the function
$\lambda f_D$ has the same properties.,
\end{reml}

\begin{lemma}
\label{lemmaL*}
Let $L \in \Pic (B)$ be a line bundle such that $c_B^* (L)=L^*$. Then,
for every divisor $D$ associated to $L$, there exists a meromorphic
function $f_D$ on $B$ such that $\Div (f_D) = D + c_B (D)$ and $f_D
= \overline{f_D \circ c_B}$.
\end{lemma}

{\bf Proof :}

By assumption, $c_B (D)$ and $-D$ are linearly equivalent. As a
consequence, there exists on $B$ a meromorphic function $f$ such that $\Div
(f) = D + c_B (D)$. Then, $g = \overline{f \circ c_B}$ is a
meromorphic function on $B$ satisfying $\Div (g) = c_B (D) + D = \Div
(f)$. Thus there exists a constant $\lambda \in \C^*$
such that
$g=\lambda f$. But, 

$$\lambda = \frac{g \circ c_B }{f \circ c_B } =
\frac{\overline{f}}{f \circ c_B } = \overline{
\Big(\frac{f}{\overline{f \circ c_B }} \Big)} =
\frac{1}{\overline{\lambda }}$$ 

Hence there exists $\theta \in \R$ such that $\lambda = \exp (2i\theta)$,
and we define $f_D = exp (i\theta) f$. $\square$

\begin{rem}
\label{remarksign}
The function $f_D$ given by lemma \ref{lemmaL*} is not unique : for
every $\lambda \in \R^*$, the function $\lambda f_D$ has the same
properties. Note that every zero or pole of 
$f_D$ on $\R B$ has even order, so that the sign of $f_D$ is constant on every
component of $\R B$.
\end{rem}

\subsection{Some particular real structures}

Let $D =
\sum_{i=1}^k n_i p_i $ be a divisor on $B$, where $p_i 
\in B$ and $n_i \in \Z$ ($i \in \{1, \dots , k\}$). We can assume
that the set $\{ p_i 
\, | \, 1 \leq i \leq k \}$ is invariant under $c_B$ (add some points
with zero coefficients to $D$ if necessary). Denote by $U_0 = B
\setminus \{ p_i \, | \, 1 \leq i \leq k \}$ and for every $i \in \{1,
\dots , k\}$, choose a holomorphic chart $(U_{p_i} , \phi_{p_i})$
such that $U_{p_i} \cap U_{p_j} = \emptyset$ if $i \neq j$, $c_B
(U_{p_i}) = U_{c_B(p_i)}$ and $\phi_{p_i} : U_{p_i} \to \Delta = \{ z \in
\C \, | \, |z| < 1 \}$ is a biholomorphism. Require in addition
that $\phi_{p_i} (p_i) 
= 0 \in \Delta$ and $\phi_{c_B(p_i)} \circ c_B \circ \phi_{p_i}^{-1} (z) =
\overline{z}$ for all $z \in \Delta$ and $i \in \{1,
\dots , k\}$ (such charts always exist, see \cite{Nat}). Such an atlas 
is called {\it compatible} with the divisor $D$ and the group $<c_B>$.
 
For every $i \in \{1, \dots , k\}$, denote by $\psi_i$ the morphism :

$$\begin{array}{rcl}
(U_{p_i} \setminus p_i) \times \C & \to & U_0 \times \C \\
(x,z) & \mapsto & (x , \phi_{p_i} (x)^{-n_i} z).
\end{array}$$
The morphisms $\psi_i$ allow to glue together the trivialisations
$U_{p_i} \times \C$, $i \in \{0, \dots , k\}$, in order to define the
line bundle $L$ associated to $D$. Such trivialisations
are called {\it compatible} with the divisor $D$ and the group $<c_B>$. 

Let $L$ (resp. $X$) be a line bundle (resp. a ruled surface) over
$B$. The real structure $c_L$ on $L$ (resp. $c_X$ on $X$) is said to
be {\it fibered over} $c_B$, or that it {\it lifts} $c_B$, if $p \circ
c_L = c_B \circ p$ (resp. $p \circ c_X = c_B \circ p$) where $p$ is
the projection $L \to B$ (resp. $X \to B$).

\begin{lemma}
\label{lemmaRL}
There exists a real structure on $L \in \Pic (B)$ which lifts $c_B$
if and only if $c_B^* (L)=L$ and for every couple $(D,f_D)$ given by lemma
\ref{lemmaL}, $f_D
\times \overline{f_D \circ c_B} = + 1$.
\end{lemma}

{\bf Proof :}

${\bf \Longrightarrow :}$ Let $s$ be a meromorphic section of $L$ and
$D = \Div (s)$. Let $c_L$ be a real structure on $L$ and $\tilde{s}=
c_L \circ s \circ c_B$. Then $\tilde{s}$ is another meromorphic
section of $L$. This implies that $\Div (\tilde{s})$ and $\Div (s)$
are linearly 
equivalent. Since $\Div (\tilde{s}) =
c_B (\Div (s))$, we deduce that $c_B^* (L) = L$. Moreover, $\tilde{s} = fs$
where $f$ is a meromorphic function on $B$ satisfying $\Div (f) = c_B (D) - D$. 
Since $s=c_L \circ \tilde{s} \circ c_B = c_L \circ (fs) \circ c_B = 
\overline{f \circ c_B} \times \tilde{s} = \overline{f \circ c_B} \times fs$,
we have $\overline{f \circ c_B} \times f = +1$. Changing the section $s$, the 
same is obtained for any couple $(D,f_D)$ given by lemma
\ref{lemmaL}.

${\bf \Longleftarrow :}$ Let $L$ be a line bundle such that $c_B^* (L)=L$ 
and $(D,f_D)$ a couple given by lemma \ref{lemmaL} such that
$f_D \times \overline{f_D \circ c_B} = + 1$. Denote
$D= \sum_{i=1}^k n_i p_i$ and let $U_0 = B
\setminus \{ p_i \, | \, 1 \leq i \leq k \}$ and $(U_{p_i} ,
\phi_{p_i})$, $i \in \{1,
\dots , k\}$, be an atlas 
compatible with the divisor $D$ and the group $<c_B>$. 

The maps 
$$\begin{array}{rcl}
U_0  \times \C & \to & U_0 \times \C \\
(x,z) & \mapsto & (c_B (x) , f_D \circ c_B (x) \overline{z}),
\end{array}$$
and for every $i \in \{1, \dots , k\}$,
$$\begin{array}{rcl}
U_{p_i}  \times \C & \to &  U_{c_B(p_i)} \times \C \\
(x,z) & \mapsto & (c_B (x) , f_D \circ c_B (x) \overline{\phi_{p_i}
(x)}^{n_{c_B(p_i)}-n_{p_i}} \overline{z}) 
\end{array}$$
glue together to form an antiholomorphic map $c_L$ on $L$. This map
lifts $c_B$ and is an involution, hence the result. $\square$

\begin{prop}
\label{propRL*}
Let $L \in \Pic (B)$ be a line bundle such that $c_B^* (L) =
L^*$. Then to every couple $(D, f_D)$ given by lemma \ref{lemmaL*} is 
associated
a real structure $c_{f_D}$ on the ruled surface $X= P(L \oplus L_0)$
which lifts $c_B$. The real part of $(X , c_{f_D})$ is orientable and
consists of $t^+$ tori, where $t^+$ is the number of components of $\R
B$ on which $f_D$ is non-negative (see remark \ref{remarksign}).
\end{prop}

\begin{rem}
\label{remarkcx+}
For the sake of simplicity, when there will not be any ambiguity on the choice
of the function $f_D$, we will denote by $c_X^+$ (resp. $c_X^-$) the
real structure $c_{f_D}$  (resp. $c_{-f_D}$).
The real part of $(X,
c_{-f_D})$ consists of $t^-$ tori, where $t^-$ is the number of components 
of $\R
B$ on which $f_D \leq 0$. Obviously, $t^+ + t^-  = \mu (\R B)$, where
$\mu (\R B)$ is the number of components of $\R B$. Thus, when $\mu (\R B)$ 
is odd,
the real structures $c_X^+$ and $c_X^-$ on $X$ cannot be
conjugated, since the numbers of components of their real parts do not
have the same parity. Nevertheless, these two real structures may
sometimes be conjugated. This situation will be studied in the next
section, proposition \ref{propconjcx+}. 
\end{rem}

{\bf Proof :}

Let $(D, f_D)$ be a couple given by 
lemma \ref{lemmaL*}, so that $f_D = \overline{f_D \circ c_B}$ and
$\Div (f_D) = D + c_B (D)$. Let $p_i 
\in B$ and $n_i \in \Z$, $i \in \{1, \dots , k\}$, be such that $D =
\sum_{i=1}^k n_i p_i $. We can assume that the set $\{ p_i
\, | \, 1 \leq i \leq k \}$ is invariant under $c_B$. Let $U_0 = B
\setminus \{ p_i \, | \, 1 \leq i \leq k \}$ and $(U_{p_i} ,
\phi_{p_i})$, $i \in \{1,
\dots , k\}$, be an atlas
compatible with the divisor $D$ and the group $<c_B>$. 

The morphisms :

$$\begin{array}{rcl}
(U_{p_i} \setminus p_i) \times \C P^1 & \to & U_0 \times \C P^1 \\
(x,(z_1 : z_0)) & \mapsto & (x , (\phi_{p_i} (x)^{-n_i} z_1 : z_0))
\end{array}$$
($i \in \{1, \dots , k\}$) allow to glue together the trivialisations
$U_{p_i} \times \C P^1$, $i \in \{0, \dots , k\}$, in order to define the
ruled surface $X$. 

Now, the maps 
$$\begin{array}{rcl}
U_0  \times \C P^1& \to & U_0 \times \C P^1 \\
(x,(z_1 : z_0)) & \mapsto & (c_B (x) , (\overline{z_0} : f_D \circ c_B
(x) \overline{z_1})), 
\end{array}$$
and for every $i \in \{1, \dots , k\}$,
$$\begin{array}{rcl}
U_{p_i}  \times \C P^1& \to &  U_{c_B(p_i)} \times \C P^1\\
(x,(z_1 : z_0)) & \mapsto & (c_B (x) , (\overline{z_0} : f_D \circ c_B
(x)\overline{\phi_{p_i} 
(x)}^{-n_{c_B(p_i)}-n_{p_i}} \overline{z_1}) 
\end{array}$$
glue together to form an antiholomorphic map $c_{f_D}$ on $X$. This map
lifts $c_B$ and is an involution. The first part of 
proposition
\ref{propRL*} is proved. 

Now, the fixed point set of $c_{f_D}$ in $U_0 \times \C
P^1$ is :
$$\{ (x, (\theta : \sqrt{f_D (x)})) \in U_0 \times \C
P^1 \, | \, x \in \R B, f_D (x) \geq 0 \text{ and } \theta \in \C ,
|\theta |=1 \}.$$
The connected components of this fixed point set are then tori or
cylinders depending on whether the corresponding component of $\R B$ is
completely included in $U_0$ or not. Similarly, the fixed point set of
$c_{f_D}$ in $U_{p_i} \times \C P^1$ is :
$$\{ (x, (\theta_i : \sqrt{f_D (x) \times x_i^{-2n_i}})) \in U_{p_i} \times \C
P^1 \, | \, x \in \R B, f_D (x) \geq 0 \text{ and } \theta_i \in \C ,
|\theta_i |=1 \},$$ 
where $x_i = \phi_{p_i} (x)$. This fixed point set is a cylinder if
$p_i \in \R B$ and is empty otherwise.

The gluing maps between these cylinders are given by $\theta =
- \theta_i$ if $x_i = \phi_{p_i} (x) < 0 $ and by $\theta =
\theta_i$ if $x_i = \phi_{p_i} (x) > 0 $. Since both $id$ and
$-id$ preserve the orientation of the circle $U^1 = \{ z \in \C \, |
\, |z| =1 \}$, the results of these gluings are always tori.
Thus, the real part of $(X, c_{f_D})$ consists only of tori and the
number of such tori is the number of components of $\R B$ on which
$f_D \geq 0$, that is $t^+$. $\square$

\section{Conjugacy classes of real structures on decomposable ruled
surfaces}

\subsection{Lifting of automorphisms of $X$}

I could not find the following proposition in the literature, so I give
it here.

\begin{prop}
\label{propaut}
Let $L$ be a complex line bundle over $B$ and $X$ be
the ruled surface $P (E)$, where $E=L \oplus L_0$. 

If $L \neq L^*$ or if $L = L_0$, then every automorphism of $X$
fibered over the identity of $B$ lifts to an automorphism of $E$.
If $L = L^*$ and $L \neq L_0$, then the automorphisms of $X$
fibered over the identity of $B$ which lift to automorphisms of $E$
form an index two subgroup of the group of automorphisms of $X$
fibered over the identity. In that case, the automorphisms of $X$
which do not lift are of the form 
$$\phi_\lambda = \left[ \begin{array}{cc}
0&\lambda s\\
s&0
\end{array} \right],$$
where $\lambda \in \C^*$ and $s$ is a non-zero meromorphic section of $L$.
\end{prop}

\begin{rem}
The automorphims $\phi_\lambda$ introduced in proposition
\ref{propaut} are holomorphic involutions of $X$.
\end{rem}

{\bf Proof :}

Denote by ${\cal O}_B^*$ the sheaf of holomorphic functions on $B$
which do not vanish and by ${\cal A}ut (E)$ (resp. ${\cal A}ut (X)$)
the sheaf of automorphisms of $E$ (resp. of $X$) fibered over the
identity. These sheafs satisfy the exact sequence :
$$1 \to {\cal O}_B^* \to {\cal A}ut (E) \to {\cal A}ut (X) \to 1$$
We deduce the following long exact sequence :
$$1 \to H^0 (B, {\cal O}_B^*) \to H^0 (B,  {\cal A}ut (E)) \to H^0 (B,
{\cal A}ut (X)) \to H^1 (B, {\cal O}_B^*) \to H^1 (B,  {\cal A}ut (E))
\to \dots$$
We are searching for the image of the morphism $H^0 (B,  {\cal A}ut (E)) \to 
H^0 (B,
{\cal A}ut (X))$. To compute this image, let us study the kernel of the
map $i_* : H^1 (B, {\cal O}_B^*) \to H^1 (B,  {\cal A}ut (E))$. 

Remember that the group $H^1 (B, {\cal O}_B^*)$ is isomorphic to $\Pic
(B)$. Such an isomorphism can be defined as follows : fix a divisor
$\sum_{j=1}^t r_j q_j$, where for $j \in \{1, \dots , t\}$, $r_j \in
\Z$ and $q_j \in  B$. Denote by $U_0 = B
\setminus \{ q_j \, | \, 1 \leq j \leq t \}$ and for every $j \in \{1,
\dots , t\}$, choose a holomorphic chart $(U_{q_j} , \phi_{q_j})$ of $B$
such that $U_{q_j} \cap U_{q_{j'}} = \emptyset$ if $j \neq j'$,
$\phi_{q_j} : U_{q_j} \to \Delta = \{ z \in
\C \, | \, |z| < 1 \}$ is a biholomorphism and $\phi_{q_j} (q_j) 
= 0 \in \Delta$. Denote by ${\cal U}$ the covering of $B$ defined by
$U_0, \dots, U_t$ and consider the following sections of ${\cal
O}_B^*$ ($j \in \{1, \dots , t\}$) :
\begin{eqnarray*}
l^1_{0j} : U_0 \cap U_j & \to & \C^* \\
x & \mapsto & \phi_{q_j} (x)^{r_j} = x_j^{r_j},
\end{eqnarray*}
where by definition $x_j = \phi_{q_j} (x) \in \Delta$. These sections
define a $1$-cocycle of $B$ with coefficient in ${\cal
O}_B^*$ and we denote with the same letter $l^1$ its
cohomology class in $H^1 ({\cal U}, {\cal O}_B^*)$ and in $H^1 (B,
{\cal O}_B^*)$. This construction defines an isomorphism between
$\Pic (B)$ and $H^1 (B, {\cal O}_B^*)$.

So let $l^1 \in H^1 (B, {\cal O}_B^*)$ be associated to the divisor
$\sum_{j=1}^t r_j q_j$. Then $m^1 = i_* (l^1)$ is the cohomology class
of the $1$-cocycle with coefficient in ${\cal A}ut (E)$ defined by the
following sections ($j \in \{1, \dots , t\}$) :
\begin{eqnarray*}
m^1_{0j} : U_0 \cap U_j & \to & {\cal A}ut (E) \\
x & \mapsto & \left[ \begin{array}{cc}
x_j^{r_j}&0\\
0&x_j^{r_j}
\end{array} \right].
\end{eqnarray*}
Suppose that $m^1 = 0 \in H^1 (B,  {\cal A}ut (E))$. Then
$\sum_{j=1}^t r_j q_j$ is of degree zero, since $0=\det (m^1) = 2l^1
\in H^1 (B, {\cal O}_B^*)$. Moreover, since the map $H^1 ({\cal
U},{\cal A}ut (E)) \to H^1 (B,  {\cal A}ut (E))$ is injective (see
\cite{Mir}, lemma $3.11$, p$294$), $m^1$
is the coboundary of a $0$-cochain given in the covering ${\cal
U}$ by the following sections ($j
\in \{0, \dots , t\}$) : 
\begin{eqnarray*}
m^0_{j} : U_j & \to & {\cal A}ut (E) \\
x & \mapsto & \left[ \begin{array}{cc}
a_j (x)&c_j (x)\\
b_j (x)&d_j (x)
\end{array} \right],
\end{eqnarray*}
where $a_j$, $d_j$ are $0$-cochains with coefficients in ${\cal O}_B$,
$c_j$ is a $0$-cochain with coefficient in ${\cal O}_B (L)$, $d_j$ is
a $0$-cochain with coefficient in ${\cal O}_B (L^*)$ and $a_j d_j -
b_jc_j$ does not vanish. Then, the equality $m^1 = \delta m^0$ can be
written :
$$\forall j \in \{1, \dots , t\}, \quad m^1_{0j}= m^0_0
(m_j^0)^{-1},$$
which rewrites as $m^0_0 = x_j^{r_j} m_j^0$ ($j \in \{1, \dots ,
t\}$).
Hence, we deduce that for $j \in \{1, \dots ,
t\}$, $a_0 = x_j^{r_j} a_j$, $d_0 = x_j^{r_j} d_j$, $b_0 = x_j^{r_j}
b_j$ and $c_0 = x_j^{r_j} c_j$. As soon as $a_0$ (resp. $d_0$) is
non-zero, this implies that $a_0$ (resp. $d_0$) is a meromorphic
function over $B$ satisfying $\Div (a_0) \geq \sum_{j=1}^t r_j q_j$
(resp. $\Div (d_0) \geq \sum_{j=1}^t r_j q_j$). Since these two
divisors are of degree zero, they are equal. So $\sum_{j=1}^t r_j q_j$
is a principal divisor and $l^1 = 0$. When $a_0 = d_0 = 0$, we deduce
that $b_0$ (resp. $c_0$) is a meromorphic
section of $L^*$ (resp. of $L$) satisfying $\Div (b_0) \geq \sum_{j=1}^t r_j q_j$
(resp. $\Div (c_0) \geq \sum_{j=1}^t r_j q_j$). Since $\deg (L) =
-\deg (L^*)$, these divisors are equal. We then deduce that $L = L^*$
and that this line bundle is associated to the divisor $\sum_{j=1}^t
r_j q_j$.

In conclusion, when $L \neq L^*$, the morphism $i_*$ is injective and
when $L = L^*$, $L \neq L_0$, the kernel of $i_*$ is included into the
subgroup of $H^1 (B, {\cal O}_B^*) = \Pic (B)$ generated by $L$, which
is of order two. In that case, it is not difficult to check that the
kernel of $i_*$ is exactly this subgroup of order two. Indeed, with
the preceding notations, it suffices to let $a_0$ and $d_0$ be equal
to $0$ and let $b_0$ and $c_0$ be equal to a same meromorphic section
of $L$. This constructs a $0$-cochain $m^0$ such that $\delta m^0 =
i_* (L)$. The first part of the proposition is proved.

To check the second part of the proposition, note that when $L = L^*
\neq L_0$, $H^0 (B , L) = H^0 (B , L^*) = 0$, so that the
automorphisms of $E = L \oplus L_0$ fibered over the identity of $B$
are of the form 
$$ \left[ \begin{array}{cc}
a&0\\
0&d
\end{array} \right],$$
where $a,d \in \C^*$. The automorphisms of $X$ fibered over the identity
which lift to $E$ are then of the form 
$$ \left[ \begin{array}{cc}
1&0\\
0&\lambda
\end{array} \right] \quad (\lambda \in \C^*).$$
It follows that the automorphisms $\phi_\lambda$ do not lift to
automorphims of $E$ and that they are the only ones with this
property. $\square$

\subsection{The conjugation's theorem}

Denote by $ c_{L_0}$ the real structure on $L_0$ defined by :
$$\begin{array}{rcl}
B \times \C & \to & B \times \C \\
(x,z) & \mapsto & (c_B (x) , \overline{z})
\end{array}$$
This real structure lifts $c_B$.

\begin{theo}
\label{theoremrealstruct}
Let $L$ be a line bundle over a smooth compact complex irreducible curve $B$
equipped with a real structure $c_B$
and let $X=P(L \oplus L_0)$ be the associated decomposable ruled
surface.

1. Suppose that $L \neq L^*$ and that there exists a real structure $c_L$ on
$L$ which lifts $c_B$. Then there exists, up to conjugation by a 
biholomorphism of $X$, one and
only one real structure on $X$ which lifts $c_B$. It is the real
structure induced by $c_L \oplus c_{L_0}$.

2. Suppose that $c_B^* (L) = L^*$. If $L \neq L^*$, then every
real structure on $X$ which 
lifts $c_B$ is conjugated to one of the two structures $c_X^+$ or
$c_X^-$ given by proposition \ref{propRL*}. The same result occurs when
$L=L_0$ or when $L = L^*$ and there is no real structure on
$L$ which lifts $c_B$.

3. Suppose that $c_B^* (L) = L = L^*$, that $L \neq L_0$ and that
there exists a real structure $c_L$ on $L$ which lifts $c_B$. Then every
real structure on $X$ which 
lifts $c_B$ is conjugated to the real structure $c_L \oplus c_{L_0}$, or
to one of the two structures $c_X^+$ or
$c_X^-$ given by proposition \ref{propRL*}.

In any other case, $X$ does not admit real structures fibered
over $c_B$.
\end{theo}

\begin{rem}
It follows from lemma \ref{lemmaRL} and remark \ref{remarkL} that when 
$\R B \neq \emptyset$,
there exists a real structure on $L$ which lifts $c_B$ if and only if
$c_B^* (L) = L$. 

Note that in the third case, the real
structures $c_X^+$ and $c_X^-$ are not conjugated to $c_L \oplus
c_{L_0}$, since they are exchanging the two disjoint holomorphic sections
of zero square of $X$ and $c_L \oplus c_{L_0}$ does not.
Note also that when $X = B \times \C P^1$, or when
$\mu (\R B)$ is odd, the real structures $c_X^+$ and $c_X^-$ on $X$
are not conjugated (see remark \ref{remarkcx+}). Nevertheless, these two
real structures may sometimes be conjugated, see
proposition \ref{propconjcx+}.
\end{rem}

\begin{prop}
\label{propexistence}
Let $L$ be a line bundle over $(B, c_B)$ and let $X=P(L \oplus L_0)$. Then
there exists a real structure on $X$ which lifts $c_B$ if and only if
there exists a real structure on $L$ which lifts $c_B$ or $c_B^* (L) = L^*$.
\end{prop}

{\bf Proof :}

${\bf \Longrightarrow :}$ To begin with, suppose that $\deg (L) \neq
0$. Then, without loss of generality, we can assume that $d= \deg (L)
>0$. The holomorphic section $e$ of $X$ associated to $L$
satisfy $e \circ e = -d <0$, since its normal bundle is $L^*$. Any
other section $\tilde{e}$ of $X$ is homologous to $e + kv$, where $k
\in \Z$ and 
$v$ is the integer homology class of a fiber. When $\tilde{e} \neq e$,
we have
$\tilde{e} \circ e \geq 0$, which means that $k \geq d$. Then
$\tilde{e} \circ \tilde{e} \geq d$ and this proves that $e$ is the only
holomorphic section  of $X$ with negative square. Thus this section is
invariant under the real structure of $X$, and so is its normal
bundle. This implies that there exists a real structure on $L^*$ which
lifts $c_B$. Using duality, there exists one on $L$ which
lifts $c_B$.

Suppose now that $\deg (L) = 0$. If $L$ is the trivial bundle, then $X
= B \times \C P^1$ and nothing has to be proved. Otherwise, the
sections of $X$ associated to $L$ and $L_0$ are the only ones with zero
squares. Indeed, a third holomorphic section with  zero
square should be disjoint from them and these three sections would
give a trivialisation of $X$. This would contradict the assumption that $X
\neq B \times \C P^1$. As a consequence, we deduce the following alternative
: either the real structure $c_X$ preserves these two sections, or it
exchanges them. In the first case, $c_X$ preserves the normal bundles
and we conclude as before. In the second case, $c_X$ exchanges the
normal bundles and so defines a morphism $\hat{c}_X : L^* \to L$,
fibered over 
$c_B$. Let $s$ be a meromorphic section of $L^*$, so that $\Div (s) =
-D$ where $D$ is a divisor associated to $L$. Then $\hat{c}_X \circ s \circ
c_B$ is a meromorphic section of $L$ and $\Div (\hat{c}_X \circ s \circ
c_B) = c_B^* (\Div (s)) = - c_B^* (D)$. Hence $c_B^* (L) = L^*$.

${\bf \Longleftarrow :}$ If there exists a real structure on $L$ which
lifts $c_B$, then taking the direct sum with $c_{L_0}$ 
we get a real structure on $L \oplus L_0$ which
lifts $c_B$. This structure induces on $X=P(L \oplus L_0)$ a real
structure which 
lifts $c_B$. If $c_B^* (L) = L^*$, the result follows from proposition
\ref{propRL*}. $\square$\\

{\bf Proof of theorem \ref{theoremrealstruct} :}

When $X = B \times \C P^1$, the second part of theorem 
\ref{theoremrealstruct} is clear. Indeed, in this case every real
structure on $X$ 
which lifts $c_B$ is the direct sum of $c_B$ and a real structure on
$\C P^1$. Moreover, the group of automorphisms of $X$ fibered over the
identity is then equal to $\{ id \} \times Aut (\C P^1)$. So the second
part of  theorem
\ref{theoremrealstruct} follows from the well known fact that, up to
conjugation, there are two real structures on $\C P^1$. Thus, from now
on, we can assume that $L \neq L_0$.
It follows from proposition \ref{propexistence} that if there exists a
real structure on $X$ which lifts $c_B$, then either there exists a
real structure $c_L$ on $L$ which
lifts $c_B$, or $c_B^* (L) = L^*$. This already proves the last line
of theorem \ref{theoremrealstruct}. We will show the theorem in three
steps.

In the first step, we will prove that if there exists a
real structure $c_L$ on $L$ which
lifts $c_B$, then every real structure on $X$ of the form $c_X \circ
\phi$, where $c_X$ is the real structure of $X$ induced by $c_L \oplus
c_{L_0}$ and $\phi$ is an automorphism of $X$ fibered over the identity
of $B$ which lifts to an automorphism of $E=L \oplus L_0$, is
conjugated to $c_X$. In the second step, we will prove that if $c_B^* (L) = L^*$,
then every real structure on $X$ of the form $c_X^+ \circ
\phi$, where $\phi$ is an automorphism of $X$ fibered over the identity
of $B$ which lifts to an automorphism of $E=L \oplus L_0$, is
conjugated either to $c_X^+$ or to $c_X^-$. Finally, in the third
step, we will prove that if $c_B^* (L) = L^* = L$, then every real
structure on $X$ of the form $c_X^+ \circ
\phi$, where $\phi$ is an automorphism of $X$ fibered over the identity
of $B$ which does not lift to an automorphism of $E=L \oplus L_0$, is
conjugated to a real structure of the form $c_L \oplus c_{L_0}$, where
$c_L$ is a real structure on $L$ which lifts $c_B$. Furthermore, this
conjugation is given by an automorphism of $X$ fibered over the identity
of $B$ which lifts to an automorphism of $E=L \oplus L_0$. In
particular, when there is no real structure on $L$ which lifts $c_B$,
every antiholomorphic map of the form $c_X^+ \circ
\phi$, where $\phi$ is an automorphism of $X$ fibered over the identity
of $B$ which does not lift to an automorphism of $E=L \oplus L_0$, is
not an involution. The theorem follows from these three steps and 
proposition \ref{propaut}.\\

{\bf First step :} Suppose that there exists a
real structure $c_L$ on $L$ which
lifts $c_B$ and let $c_X$ be the real structure of $X$ induced by $c_L \oplus
c_{L_0}$. Let $\tilde{c}_X$ be another real structure on $X$ which is  of the 
form $c_X \circ
\phi$, where $\phi$ is an automorphism of $X$ fibered over the identity
of $B$ which lifts to an automorphism of $E=L \oplus L_0$. The aim of
this first step  is to prove that $c_X$ and  $\tilde{c}_X$ are
conjugated.

Let $\Phi$ be an automorphism of $E= L \oplus L_0$ which
lifts $\phi$. Then $\Phi \in End (E) = E \otimes E^* = L \oplus L^* \oplus L_0
\oplus L_0$. Thus there exist $a,d \in \C^*$, $b \in H^0 (B , L^*)$
and $c \in H^0 (B , L)$ such that 
$$\Phi = \left[ \begin{array}{cc}
a&c\\
b&d
\end{array} \right]$$
By assumption, the line bundle $L$ is not trivial, so that either $L$ or
$L^*$ has no non-zero holomorphic section. Without loss of generality,
we can assume that it is $L$, so that $c=0$ and 
$$\Phi = \left[ \begin{array}{cc}
a&0\\
b&d
\end{array} \right].$$
By assumption, $\tilde{c}_X^2 = id$, which implies that $c_X \circ
\phi \circ c_X 
= \phi^{-1}$. So there exists $\lambda \in \C^*$ such that $c_E \circ
\Phi \circ c_E = \lambda \Phi^{-1}$. But 
$$\Phi^{-1} = \frac{1}{ad}\left[ \begin{array}{cc}
d&0\\
-b&a
\end{array} \right],$$
and
$$c_E \circ \Phi \circ c_E = \left[ \begin{array}{cc}
\overline{a}&0\\
c_{L_0} \circ b \circ c_L&\overline{d}
\end{array} \right].$$
Put $\tilde{\lambda} = \frac{1}{ad} \lambda$, we have $\tilde{\lambda}
d=\overline{a}$, $\tilde{\lambda} a=\overline{d}$ and
$-\tilde{\lambda} b= c_{L_0} \circ b \circ c_L$. The two first conditions
imply that $|\tilde{\lambda}| = 1$. Thus there exists $\theta \in \R$
such that $\tilde{\lambda} = \exp (2i \theta )$. So the previous
conditions can be rewritten as $\exp (i \theta ) d = \overline{\exp (i
\theta )a} $, $\exp (i \theta ) a= \overline{\exp (i
\theta )d} $ and $- \exp (i \theta )b= c_{L_0} \circ (\exp (i \theta )b) 
\circ c_L$.
Hence we can assume that 
$$\Phi = \left[ \begin{array}{cc}
a&0\\
b&d
\end{array} \right],$$
where $d=\overline{a}$ and $b=-c_{L_0} \circ b \circ c_L$ (replace $\Phi$
by $\exp (i \theta ) \Phi$ which also lifts $\phi$).

Now, denote by $\Psi$ the automorphism of $E$ defined by 
$$\Psi = \left[ \begin{array}{cc}
1&0\\
\frac{1}{2}b&\overline{a}
\end{array} \right]. $$
Then 
$$\Psi^{-1} = \frac{1}{\overline{a}} \left[ \begin{array}{cc}
\overline{a}&0\\
-\frac{1}{2}b  &1
\end{array} \right] \text{, and}$$
\begin{eqnarray*}
\Psi^{-1} \circ c_E \circ \Psi &=& \frac{1}{\overline{a}}\left[
\begin{array}{cc} 
\overline{a}c_L&0\\
-\frac{1}{2}b\circ c_L + \frac{1}{2}c_{L_0} \circ {b}  &a c_{L_0}
\end{array} \right] \\
&=&\frac{1}{\overline{a}}\left[
\begin{array}{cc} 
\overline{a}c_L&0\\
c_{L_0} \circ {b}  &a c_{L_0}
\end{array} \right] \quad \text{since } -b \circ c_L = c_{L_0} \circ {b}\\
&=&\frac{1}{\overline{a}} c_E \circ \Phi.
\end{eqnarray*}
Denote by $\psi$ the automorphism of $X$ induced by $\Psi$, we then
deduce that $\psi^{-1} \circ c_X \circ \psi = \tilde{c}_X$, which was
the aim of this first step.\\

{\bf Second step :} Suppose that $c_B^* (L) = L^*$ and fix a real
structure $c_X^+$ on $X$ given by proposition \ref{propRL*} (see remark 
\ref{remarkcx+}). Let
$\tilde{c}_X$ be another real structure on $X$ which is  of the form $c_X^+ \circ
\phi$, where $\phi$ is an automorphism of $X$ fibered over the identity
of $B$ which lifts to an automorphism of $E=L \oplus L_0$. The aim of
this second step is to prove that $\tilde{c}_X$  is
conjugated either to $c_X^+$ or to $c_X^-$.
Let $\Phi$ be an automorphism of $E= L \oplus L_0$ which
lifts $\phi$. Since $\deg (L) = 0$ and $L$ is not trivial, $H^0 (B,L) = 
H^0 (B ,
L^*) = 0$. As a consequence, there exists $a,d \in \C^*$ such that 
$$\Phi = \left[ \begin{array}{cc}
a&0\\
0&d
\end{array} \right]. $$
Since $\tilde{c}_X^2 = id$, $\frac{a}{d} \in \R^*$
and we can assume that $a=1$, $d \in \R^*$ (replace $\Phi$ by
$\frac{1}{a} \Phi$). Let $\psi$ be the automorphism of $X$ defined by 
$$\psi = \left[ \begin{array}{cc}
1&0\\
0&\delta
\end{array} \right], $$
where $\delta = \frac{1}{\sqrt{|d|}}$. Then $\psi$ conjugates
$\tilde{c}_X$ to one of the two real structures ${c}_X^+$ or
${c}_X^-$. \\

{\bf Third step :} Suppose that $c_B^* (L) = L^*$ and fix a real
structure $c_X^+$ on $X$ given by proposition \ref{propRL*}. Let
$\tilde{c}_X$ be another real structure on $X$ which is  of the form 
$c_X^+ \circ
\phi$, where $\phi$ is an automorphism of $X$ fibered over the identity
of $B$ which does not lift to an automorphism of $E=L \oplus L_0$. The aim of
this third step is to prove that $\tilde{c}_X$  is
conjugated to a real structure of
the form $c_L \oplus c_{L_0}$ where $c_L$ is a real structure on $L$
which lifts $c_B$. 
 Note that the automorphism $\phi$ and the involution $c_X^+$ both
exchange the sections of $X$ associated to $L$ and $L_0$. Thus
$\tilde{c}_X$ preserves these two sections. As a consequence, it
preserves also the normal bundles of these sections and so induce a
real structure on the line bundle $L$ which lifts $c_B$. Consider then
the real structure $c_L \oplus c_{L_0}$ on $X$, it follows from the
first and the second step that it is conjugated to $\tilde{c}_X$ by an
automorphism of $X$ which lifts to an automorphism of $E$. $\square$

\subsection{When are $c_X^+$ and $c_X^-$ conjugated ?}

In this subsection is given a sufficient condition for $c_X^+$ and $c_X^-$ 
to be
conjugated (see proposition \ref{propconjcx+}). One important example
where this occurs is given by corollary \ref{corspin}.

\begin{prop}
\label{propconjcx+}
Let $L$ be a line bundle over $(B, c_B)$ such that $c_B^* (L) = L^*$
and let $X=P(L \oplus L_0)$ be the associated ruled
surface. Let $(D, f_D)$ be a couple given by lemma \ref{lemmaL*}
and $c_{f_D}$, $c_{-f_D}$ be the associated real structures of $X$
(see proposition \ref{propRL*}). Suppose that
there exists $\varphi \in Aut (B)$ of finite order such that 
$\varphi \circ c_B = c_B \circ \varphi$ and :

a. either $\varphi^* (L) = L$ and there exists a meromorphic function
$g$ on $B$ such that $\Div (g) = \varphi (D) - D$ and $f_D \circ
\varphi \times g \circ \varphi \times \overline{g \circ c_B \circ
\varphi} = -f_D$,

b. or $\varphi^* (L) = L^*$ and there exists a meromorphic function
$h$ on $B$ such that $\Div (h) = \varphi (D) + D$ and $h \circ \varphi
\times \overline{h \circ c_B \circ 
\varphi} = -f_D \times f_D \circ \varphi$.

Then, the real structures $c_{f_D}$ and $c_{-f_D}$ are conjugated in $X$.
\end{prop}

\begin{rem}
When $\R B \neq \emptyset$, the conditions $a.$ and $b.$ can be
replaced by $\varphi^* (L) \in \{ L , L^* \}$ and there exists $x \in
\R B$ such that $f_D \times f_D \circ \varphi (x) <0$. Indeed, it is
not difficult to check that in the situation $a$, there always exists
a meromorphic function
$g$ on $B$ such that $\Div (g) = \varphi (D) - D$ and $f_D \circ
\varphi \times g \circ \varphi \times \overline{g \circ c_B \circ
\varphi} = \epsilon f_D$ where $\epsilon = \pm 1$. Similarly, in the
situation $b$, there always exists 
a meromorphic function $h$ on $B$ such that $\Div (h) = \varphi (D)
+ D$ and $h \circ \varphi 
\times \overline{h \circ c_B \circ  
\varphi} = \epsilon f_D \times f_D \circ \varphi$, where $\epsilon =
\pm 1$.
Hence, conditions $a$ or $b$ are equivalent to require that $\epsilon
= -1$, which is equivalent, when $\R B \neq \emptyset$, to require
that there exists $x \in
\R B$ such that $f_D \times f_D \circ \varphi (x) <0$.

Note that when $g(B) \geq 2$, the conditions given by proposition
 \ref{propconjcx+} 
are in fact necessary and sufficient for $c_{f_D}$ and $c_{-f_D}$ 
to be conjugated, but this will not be needed in what follows.
\end{rem}

\begin{cor}
\label{corspin}
Let $g \geq 1$ be an odd integer. Then there exists a smooth compact
irreducible real algebraic curve $(B,c_B)$ of genus $g$ and empty real
part together with a complex line bundle $L$ over $B$ satisfying
$c_B^* (L) = L^*$, such that the real structures $c_{X}^+$ and $c_{X}^-$
on $X = P(L \oplus L_0)$ are conjugated.
\end{cor}

\newpage
{\bf Proof :}

Let us consider first the case $g=1$. Let $B$ be the elliptic curve $\C
/ \Z[i]$ equipped with the real structure $c_B (z) = \overline{z} +
\frac{1}{2}$, so that $\R B = \emptyset$.
Let $p_0 = 0$, $q_0 = \frac{1}{2}$, $p_1 = \frac{i}{2}$ and $q_1
=\frac{1}{2} +  \frac{i}{2}$.
$$\vcenter
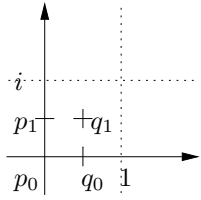
{\hbox{\begin{picture}(0,0)%
\epsfig{file=realstr1.pstex}%
\end{picture}%
\setlength{\unitlength}{3158sp}%
\begingroup\makeatletter\ifx\SetFigFont\undefined%
\gdef\SetFigFont#1#2#3#4#5{%
  \reset@font\fontsize{#1}{#2pt}%
  \fontfamily{#3}\fontseries{#4}\fontshape{#5}%
  \selectfont}%
\fi\endgroup%
\begin{picture}(1524,1524)(2389,-4573)
\put(2476,-4486){\makebox(0,0)[lb]{\smash{\SetFigFont{10}{12.0}{\rmdefault}{\mddefault}{\updefault}$p_0$}}}
\put(3076,-4036){\makebox(0,0)[lb]{\smash{\SetFigFont{10}{12.0}{\rmdefault}{\mddefault}{\updefault}$q_1$}}}
\put(2476,-3736){\makebox(0,0)[lb]{\smash{\SetFigFont{10}{12.0}{\rmdefault}{\mddefault}{\updefault}$i$}}}
\put(3301,-4486){\makebox(0,0)[lb]{\smash{\SetFigFont{10}{12.0}{\rmdefault}{\mddefault}{\updefault}$1$}}}
\put(3001,-4486){\makebox(0,0)[lb]{\smash{\SetFigFont{10}{12.0}{\rmdefault}{\mddefault}{\updefault}$q_0$}}}
\put(2476,-4036){\makebox(0,0)[lb]{\smash{\SetFigFont{10}{12.0}{\rmdefault}{\mddefault}{\updefault}$p_1$}}}
\end{picture}
}}$$
Let $D = p_1 - p_0$ and denote by $L$
the associated complex line bundle over $B$. Then $c_B^* (L)=L =
L^*$. Denote by $\varphi$ the involutive
automorphism of $B$ defined by $\varphi (z) = z + \frac{1}{2}$. Then
$\varphi \circ c_B = c_B \circ \varphi$ and $\varphi^* (L) = L$. We will 
prove that
$\varphi$ satisfies condition $a$ of proposition \ref{propconjcx+}.

For this, let $f$ be a meromorphic function on $B$ given by lemma
\ref{lemmaL*}, such that $\overline{f \circ c_B} = f$ and $\Div (f) =
D + c_B (D)$. Then $ f \circ \varphi = f$. Indeed,
there exists a holomorphic section $s$ of the line bundle $L$ such
that $\Div (s) = D$ and $s \otimes (s \circ \varphi) = f$. Thus $f \circ
\varphi = (s \circ \varphi) \otimes s = s \otimes (s \circ \varphi) = f$.
Now let
$g$  be a meromorphic function on $B$ such that $\Div
(g) = \varphi (D) - D = q_1 - p_1 - q_0 + p_0$ and $g \times 
\overline{g
\circ c_B} = -1$. Such a function is given by lemma \ref{lemmaL} and
\cite{GH}, proposition 
$2.2$, since $D$ belong to the nontrivial component of the 
real
part of $(\Jac (B) , c_B)$. Then $f \circ \varphi \times
g \circ \varphi \times \overline{g \circ c_B \circ \varphi} = -f$, so
that the condition $a$ of proposition \ref{propconjcx+} is
satisfied. We deduce that the real structures $c_{X}^+$ and
$c_{X}^-$ on $X = P(L \oplus L_0)$ defined by $f$ and $-f$ 
(see proposition \ref{propRL*}) are conjugated.

Now, let us consider the case $g=2k+1$, $k \geq 1$. For $j \in \{ 0, \dots , 
2k-1 \}$, denote by $\tilde{p}_j = \frac{j}{2k} i \in B$ and $\tilde{q}_j = 
\frac{1}{2} + \frac{j}{2k} i \in B$ (so that $p_1 = \tilde{p}_k$ and 
$q_1 = \tilde{q}_k$). Denote by $B_k$ the double 
covering of $B$ ramified over the $4k$ points $\tilde{p}_j$, $\tilde{q}_j$,
$j \in \{0 , \dots, 2k-1 \}$. This covering can be chosen so that its 
characteristic class in
$H^1 (B \setminus \{\tilde{p}_j,\tilde{q}_j \, | \, j \in \{0 , \dots, 2k-1 
\}  \} ; \Z / 2\Z)$ is Poincar\'e dual to the sum of the $2k$ segments
$\{ (0,t) \, | \, t \in ]\frac{2j}{2k} , \frac{2j+1}{2k}[,
j \in \{0 , \dots, k-1 \} \}$ and $\{ (\frac{1}{2},t)
\, | \, t \in ]\frac{2j}{2k} , \frac{2j+1}{2k}[,
j \in \{0 , \dots, k-1 \} \}$. 
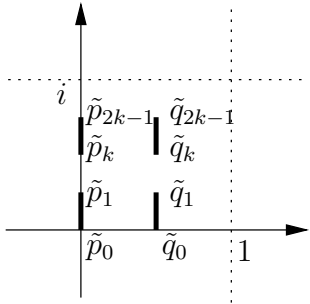
$$\vcenter{\hbox{\begin{picture}(0,0)%
\epsfig{file=realstr3.pstex}%
\end{picture}%
\setlength{\unitlength}{4144sp}%
\begingroup\makeatletter\ifx\SetFigFont\undefined%
\gdef\SetFigFont#1#2#3#4#5{%
  \reset@font\fontsize{#1}{#2pt}%
  \fontfamily{#3}\fontseries{#4}\fontshape{#5}%
  \selectfont}%
\fi\endgroup%
\begin{picture}(1824,1824)(3139,-7948)
\put(3466,-6721){\makebox(0,0)[lb]{\smash{\SetFigFont{12}{14.4}{\rmdefault}{\mddefault}{\updefault}$i$}}}
\put(4546,-7666){\makebox(0,0)[lb]{\smash{\SetFigFont{12}{14.4}{\rmdefault}{\mddefault}{\updefault}$1$}}}
\put(3646,-7621){\makebox(0,0)[lb]{\smash{\SetFigFont{12}{14.4}{\rmdefault}{\mddefault}{\updefault}$\tilde{p}_0$}}}
\put(4096,-7621){\makebox(0,0)[lb]{\smash{\SetFigFont{12}{14.4}{\rmdefault}{\mddefault}{\updefault}$\tilde{q}_0$}}}
\put(3646,-7036){\makebox(0,0)[lb]{\smash{\SetFigFont{12}{14.4}{\rmdefault}{\mddefault}{\updefault}$\tilde{p}_k$}}}
\put(3646,-7306){\makebox(0,0)[lb]{\smash{\SetFigFont{12}{14.4}{\rmdefault}{\mddefault}{\updefault}$\tilde{p}_1$}}}
\put(3646,-6811){\makebox(0,0)[lb]{\smash{\SetFigFont{12}{14.4}{\rmdefault}{\mddefault}{\updefault}$\tilde{p}_{2k-1}$}}}
\put(4141,-6811){\makebox(0,0)[lb]{\smash{\SetFigFont{12}{14.4}{\rmdefault}{\mddefault}{\updefault}$\tilde{q}_{2k-1}$}}}
\put(4141,-7036){\makebox(0,0)[lb]{\smash{\SetFigFont{12}{14.4}{\rmdefault}{\mddefault}{\updefault}$\tilde{q}_k$}}}
\put(4141,-7306){\makebox(0,0)[lb]{\smash{\SetFigFont{12}{14.4}{\rmdefault}{\mddefault}{\updefault}$\tilde{q}_1$}}}
\end{picture}
}}$$
Denote by $\pi_k : B_k \to B$ the
projection associated to the covering. The automorphism $\varphi $ of $B$
lifts to an automorphism $\varphi_k$ of $B_k$ such that $\varphi \circ \pi_k
= \pi_k \circ \varphi_k$. Similarly, the real structure $c_B$ lifts to a
real structure $c_{B_k}$ on $B_k$ such that $c_B \circ \pi_k = \pi_k
\circ c_{B_k}$ and $\R B_k = \emptyset$. Denote by $L_k = \pi_k^*
(L)$. This bundle satisfies $c_{B_k}^* (L_k) = L_k = L_k^* =
\varphi_k^* (L_k)$. Finally, denote by $f_k = f \circ \pi_k$
and $g_k = g \circ \pi_k$. Then $f_k = \overline{f_k
\circ c_{B_k}}$ and $g_k \times \overline{g_k
\circ c_{B_k}} = -1$. Moreover, $\Div (f_k) =c_{B_k}^* (D_k) +
D_k$, where $D_k = \pi_k^* (D) = 2p_1 - 2p_0$, and $\Div (g_k) =
\varphi_k^* (D_k) - D_k$. We have, $f_k \circ \varphi_k \times
g_k \circ \varphi_k \times \overline{g_k \circ c_{B_k} \circ \varphi_k} 
= -f_k$, so
that the condition $a$ of proposition \ref{propconjcx+} is
satisfied. We deduce that the real structure $c_{X_k}^+$ and
$c_{X_k}^-$ on $X_k = P(L_k \oplus L_0)$ defined by $f_k$ and $-f_k$ 
(see proposition \ref{propRL*}) are conjugated. $\square$\\

{\bf Proof of proposition \ref{propconjcx+} :}

Denote $D =
\sum_{i=1}^k n_i p_i $, where $p_i 
\in B$ and $n_i \in \Z$, $i \in \{1, \dots , k\}$. We can assume 
that the set $\{ p_i 
\, | \, 1 \leq i \leq k \}$ is invariant under $\varphi$ (add some points
with zero coefficients to $D$ if necessary). Denote by $U_0 = B
\setminus \{ p_i \, | \, 1 \leq i \leq k \}$ and for every $i \in \{1,
\dots , k\}$, choose some holomorphic chart $(U_{p_i} , \phi_{p_i})$
such that $U_{p_i} \cap U_{p_j} = \emptyset$ if $i \neq j$, $\varphi
(U_{p_i}) = U_{\varphi(p_i)}$ and $\phi_{p_i} : U_{p_i} \to \Delta = \{ z \in
\C \, | \, |z| < 1 \}$ is a biholomorphism. We require in addition
that $\phi_{p_i} (p_i) 
= 0$ and 
$$\begin{array}{rcl}
\phi_{\varphi(p_i)} \circ \varphi \circ \phi_{p_i}^{-1} : \Delta & \to &
\Delta\\
x & \mapsto & 
\exp(\frac{2i\pi}{m_i}) x \text{ if $ p_i$ is
a fixed point of order $m_i$ of $\varphi$.} 
\end{array}$$
(we put $m_i = 1$ if 
$\varphi (p_i) \neq p_i$. This atlas and these trivialisations are
compatible with $D$ and the group $<\varphi>$. It always exists, see
\cite{Nat}.) 

For every $i \in \{1, \dots , k\}$, denote by $\psi_i$ the morphism :

$$\begin{array}{rcl}
(U_{p_i} \setminus p_i) \times \C P^1 & \to & U_0 \times \C P^1\\
(x,(z_1 : z_0)) & \mapsto & (x , (\phi_{p_i} (x)^{-n_i} z_1 : z_0)).
\end{array}$$
The morphisms $\psi_i$ allow to glue together the trivialisations
$U_{p_i} \times \C P^1$, $i \in \{0, \dots , k\}$, in order to define the
ruled surface $X$. 

Now suppose we are in the case a. Let $g$ be the meromorphic function on 
$B$ such that
$\Div (g) = \varphi (D) - D$ and $f_D \circ 
\varphi \times g \circ \varphi \times \overline{g \circ c_B \circ
\varphi} = -f_D$. Consider the maps : 
$$\begin{array}{rcl}
U_0  \times \C P^1 & \to & U_0 \times \C P^1\\
(x,(z_1 : z_0)) & \mapsto & 
(\varphi (x), (g \circ \varphi (x) z_1 : z_0)), 
\end{array}$$
and for every $i \in \{1, \dots , k\}$,
$$\begin{array}{rcl}
U_{p_i}  \times \C P^1& \to &  U_{p_j} \times \C P^1\\
(x,(z_1 : z_0)) & \mapsto & 
(\varphi (x), (g \circ \varphi (x) \phi_{p_i} (x)^{n_j - n_i}
\exp(\frac{2i\pi }{m_i})z_1 : z_0)),
\end{array}$$
where $p_j$ denotes the point $\varphi (p_i)$. These maps
glue together to form an element $\Phi_g \in Aut (X)$ fibered over
$\varphi$.

The map $\Phi_g^{-1}$ is given by :
$$\begin{array}{rcl}
U_0  \times \C P^1 & \to & U_0 \times \C P^1\\
(x,(z_1 : z_0)) & \mapsto &  (\varphi^{-1} (x), (z_1 : g(x) z_0)). 
\end{array}$$
And the map $c_X^-$ is given by :
$$\begin{array}{rcl}
U_0  \times \C P^1 & \to & U_0 \times \C P^1\\
(x,(z_1 : z_0)) & \mapsto &  (c_B (x), (\overline{z_0} : -f_D \circ
c_B (x) \overline{z_1})). 
\end{array}$$
Thus $\Phi_g^{-1} \circ c_X^- \circ \Phi_g$ is given in this
trivialisation by :
$$\begin{array}{rcl}
U_0  \times \C P^1 & \to & U_0 \times \C P^1\\
(x,(z_1 : z_0)) & \mapsto &  (c_B (x), (\overline{z_0} : -f_D \circ
c_B \circ \varphi (x) \times \overline{g \circ \varphi (x)} \times g \circ
c_B \circ \varphi (x) \overline{z_1})). 
\end{array}$$
Since $f_D \circ 
\varphi \times g \circ \varphi \times \overline{g \circ c_B \circ
\varphi} = -f_D$, we conclude that $\Phi_g^{-1} \circ c_X^- \circ \Phi_g = 
c_X^+$.

Suppose now we are in the case b. Let $h$ be the meromorphic function on $B$ 
such that
$\Div (h) = \varphi (D) + D$ and $h \circ \varphi
\times \overline{h \circ c_B \circ 
\varphi} = -f_D \times f_D \circ \varphi$.
Consider then the maps : 
$$\begin{array}{rcl}
U_0  \times \C P^1 & \to & U_0 \times \C P^1\\
(x,(z_1 : z_0)) & \mapsto & 
(\varphi (x), (z_0 : h \circ \varphi (x) z_1)) 
\end{array}$$
and for all $i \in \{1, \dots , k\}$,
$$\begin{array}{rcl}
U_{p_i}  \times \C P^1& \to &  U_{p_j} \times \C P^1\\
(x,(z_1 : z_0)) & \mapsto & 
(\varphi (x), (z_0 : h \circ \varphi (x) \phi_{p_i} (x)^{-n_i - n_j}
\exp(-\frac{2i\pi }{m_i})z_1)),
\end{array}$$
where $p_j$ denotes the point $\varphi (p_i)$. These maps
glue together to form an element $\Phi_h \in Aut (X)$ fibered over
$\varphi$.

The map $\Phi_h^{-1}$ is given by :
$$\begin{array}{rcl}
U_0  \times \C P^1 & \to & U_0 \times \C P^1\\
(x,(z_1 : z_0)) & \mapsto &  (\varphi^{-1} (x), (z_0 : h(x) z_1)). 
\end{array}$$
And the map $c_X^-$ is given by :
$$\begin{array}{rcl}
U_0  \times \C P^1 & \to & U_0 \times \C P^1\\
(x,(z_1 : z_0)) & \mapsto &  (c_B (x), (\overline{z_0} : -f_D \circ
c_B (x) \overline{z_1})). 
\end{array}$$
Thus $\Phi_h^{-1} \circ c_X^- \circ \Phi_h$ is given in this
trivialisation by :
$$\begin{array}{rcl}
U_0  \times \C P^1 & \to & U_0 \times \C P^1\\
(x,(z_1 : z_0)) & \mapsto &  (c_B (x), ( -f_D \circ
c_B \circ \varphi (x) \overline{z_0} : \overline{h \circ \varphi (x)}
\times h \circ 
c_B \circ \varphi (x) \overline{z_1})). 
\end{array}$$
Since $h \circ \varphi
\times \overline{h \circ c_B \circ 
\varphi} = -f_D \times f_D \circ \varphi$, we conclude that 
$\Phi_h^{-1} \circ c_X^- \circ \Phi_h =  c_X^+$. $\square$

\section{Deformation classes of real structures on ruled surfaces}

\subsection{The real part of $(\Jac (B) , -c_B^*)$}
\label{subsectionpartition}

Remember the following well known result (see \cite{GH}, propositions
3.2 and 3.3 for instance) : 

\begin{prop}
\label{propjac}
Let $(B , c_B)$ be a smooth compact irreducible real algebraic curve. The 
Jacobian
$\Jac (B)$ of $B$ is equipped with the real structure $-c_B^*$. Then
if $\R B \neq \emptyset$, the real part of $(\Jac (B) ,
-c_B^*)$ has $2^{\mu (\R B) - 1}$ connected components, where $\mu (\R
B)$ is the number of components of $\R B$. 
If $\R B = \emptyset$, the real part of $(\Jac (B) , -c_B^*)$ is
connected if $g(B)$ is even and consists of two connected components 
otherwise. $\square$
\end{prop}
(Note that multiplication of $c_B^*$ by $-1$ does not
change the topology of the real part of $\Jac (B)$.)\\

Let $L$ be a complex line bundle over $B$ such that $c_B^* (L) = L^*$,
that is an element of the real part of $(\Jac (B) , -c_B^*)$, where
$\Jac (B)$ is identified with the part of $\Pic (B)$ of degree zero.
Let $(D, f_D)$ be a couple
given by lemma \ref{lemmaL*}. The function $f_D$ is real and
of constant sign on every component of $\R B$, thus it induces a
partition of $\R B$ in two elements $\R B \cap \overline{f_D^{-1} (\R_+^* )}$
and $\R B \cap \overline{f_D^{-1} (\R_-^*)}$. It follows from theorem
\ref{theoremrealstruct} that this partition only depends on the bundle
$L$ and not on the choice of $(D, f_D)$, since it corresponds to the
projections on $\R B$ of the real parts of $(P (L \oplus L_0) ,
c_X^+)$ and $(P (L \oplus L_0) , c_X^-)$. For the same reason, this
partition actually only depends on the connected component of the real
part of $(\Jac (B) , -c_B^*)$ and hence is an invariant
associated to these components. Note that when $\R B \neq \emptyset$
has $\mu (\R
B)$ components, the number of partitions of $\R B$ in two elements is
$2^{\mu (\R B) - 1}$. 

\begin{lemma}
\label{lemmapartition}
When $\R B \neq \emptyset$, the partitions associated to the real components
of $(\Jac (B) , -c_B^*)$ establish 
a bijection between the set of these components and the set of partitions of 
$\R B$ in two elements.
\end{lemma}

{\bf Proof :}

Let $L$ and $L'$ be two complex line bundles which
belong to $\R \Jac (B)$ and such that their associated partitions of
$\R B$ are the same. We will prove that they belong to the same
component of $\R \Jac (B)$. The result follows, since the
``partition'' map is then injective and hence bijective for
cardinality reasons.

Let $D$ (resp. $D'$) be a divisor associated to $L$ (resp. $L'$). Let
$f_D$ (resp. $f_{D'}$) be a non-zero meromorphic function on $B$ such
that $\overline{f_D \circ c_B} = f_D$ (resp. $\overline{f_{D'} \circ
c_B} = f_{D'}$) and $\Div (f_D) =  D + c_B (D)$ (resp. $\Div (f_{D'})
=  D' + c_B (D')$). It follows from lemma \ref{lemmaL*} that such
meromorphic functions exist. Since the partitions of $L$ and $L'$ are
the same, we can assume that $f_D$ and $f_{D'}$ have the same signs on
every components of $\R B$ (replace $f_{D'}$ by $-f_{D'}$
otherwise). For every $t \in [0,1]$, let $g_t = (1-t)f_D + t
f_{D'}$. Then $g_0 = f_D$, $g_1 = f_{D'}$ and for every $t \in [0,1]$,
$\overline{g_t \circ c_B} = g_t$. Moreover, for every $t \in [0,1]$,
$g_t$ is non-zero and of constant sign on each component of $\R B$. Thus 
every real zero
and real pole of $g_t$ is of even order. This implies that there
exists a continuous path $(D_t)_{t \in [0,1]}$ of divisors such that
$D_0 = D$ and for every $t \in [0,1]$, $\Div (g_t) = D_t + c_B
(D_t)$. In particular, $L$ and $L_1$ are in the same component of $\R
\Jac (B)$, where $L_1$ is the complex line bundle associated to
$D_1$. It suffices then to prove that $L_1$ and $L'$ lie in the same 
component of $\R
\Jac (B)$.

But $D_1 + c_B (D_1) = D' + c_B (D') = \Div (g_1)$. So the divisor
$E=D_1 - D'$ satisfy $c_B (E) = -E$. Thus there exist $k \in \N$ and
$p^1, \dots , p^k \in B$ such that $E = \sum_{i=1}^k n_i (p^i - c_B
(p^i))$. For every $i \in \{ 1, \dots , k \}$, choose a continuous
path $(p^i_\tau)_{\tau \in [0,1]}$ such that $p^i_0 = p^i$ and $p^i_1
\in \R B$. For every $\tau \in [0,1]$, let $E_\tau = \sum_{i=1}^k n_i 
(p^i_\tau - c_B
(p^i_\tau))$. Then $E_0 = E$, $E_1 = 0$ and for every  $\tau \in
[0,1]$, $c_B (E_\tau) = -E_\tau$. The path $F_\tau = D' + E_\tau$ is a
continuous path of divisors such that $F_0 = D_1$, $F_1 = D'$ and for every  
$\tau \in
[0,1]$, $F_\tau + c_B (F_\tau) = \Div (g_1)$. This implies that the
bundles $L_1$ and $L'$ belong to the same component of $\R
\Jac (B)$, hence the result. $\square$

\subsection{The topological type of a real ruled surface}
\label{subsectiontoptype}

Remember that to every smooth compact irreducible real algebraic
curve $(B,c_B)$ is associated a triple $(g, \mu , \epsilon)$, called the
{\it topological type} of $(B,c_B)$, where $g$ is the genus of $B$,
$\mu$ is
the number of connected components of $\R B$ and $\epsilon = 1$
(resp. $\epsilon = 0$) if $B$ is dividing (resp. if $B$ is
non-dividing). Two smooth compact irreducible real algebraic
curves are in the same deformation class if and only if they have the
same topological type (see \cite{Nat2}). Moreover, there exists a
smooth compact irreducible real algebraic 
curve of topological type $(g, \mu , \epsilon)$ if and only if
$\epsilon = 0$ and $0 \leq \mu \leq g$ or $\epsilon = 1$, $1 \leq \mu
\leq g+1$ and $\mu = g+1 \mod (2)$.

Except from the ellipsoid, that is $\C P^1 \times \C P^1$ equipped
with the real structure $(x,y) \mapsto (\overline{y} , \overline{x})$,
for every real structure $c_X$ on a ruled surface $p : X \to B$, there
exists a real structure $c_B$ on the base $B$ such that $p \circ c_X =
c_B \circ p$. In particular, the connected components of $\R X$ are
tori or Klein bottles. Note also that in the case of $\C P^1 \times \C P^1$,
the ruling given by the projection $p$ is not unique, whereas it is for
any other ruled surface. Since real structures on rational ruled surfaces are 
well
known (see theorem \ref{theodefrat}), we will assume from now on that {\bf 
the genus
of the base is non-zero}. So let $(X, c_X)$ be a real non-rational ruled 
surface of
base $(B,c_B)$. The {\it topological
type} of $(X, c_X)$ is by definition the quintuple $(t,k,g, \mu ,
\epsilon)$, where $(g, \mu , \epsilon)$ is the topological
type of $(B,c_B)$, $k$ is the number of Klein bottles of $\R X$ and $t$ the
number of tori of $\R X$. Obviously $t,k \geq 0$ and $t+k \leq \mu$. A
quintuple $(t,k,g, \mu ,\epsilon)$ is called {\it allowable} if $t,k
\geq 0$, $t+k \leq \mu$, $g \geq 1$ and either $\epsilon = 0$ and $0 \leq \mu \leq
g$ or $\epsilon = 1$, $1 \leq \mu 
\leq g+1$ and $\mu = g+1 \mod (2)$.

\begin{prop}
\label{proposexistence}
There exists a real ruled surface of topological type $(t,k,g, \mu
,\epsilon)$ if and only if the quintuple $(t,k,g, \mu ,\epsilon)$ is allowable.
\end{prop}

{\bf Proof :}

If $(t,k,g, \mu ,\epsilon)$ is the topological type of a real ruled
surface, then the quintuple $(t,k,g, \mu ,\epsilon)$ is clearly
allowable. Now, 
let $(t,k,g, \mu ,\epsilon)$ be an  allowable quintuple. It is well
known (see \cite{Nat2} for instance) that there exists a smooth
compact connected real algebraic 
curve $(B,c_B)$ whose topological type is $(g, \mu ,
\epsilon)$. If $\mu = 0$, the ruled surface $(B \times \C P^1 , c_B 
\times conj)$,
where $conj$ is a real structure on $\C P^1$, is of topological type 
$(0,0,g,0,0)$. 
If $\mu \neq 0$, choose a partition ${\cal P}$ of $\R B$ in two elements
 such that one
of them contains $t+k$ components of $\R B$ and the other one $\mu - t
-k$. It follows from lemma \ref{lemmapartition} that there exists a
line bundle $L$ over $B$ such that $c_B^* (L) = L^*$ and the partition
associated to $L$ is ${\cal P}$. Thus, it follows from proposition
\ref{propRL*} that there exists a real structure $c_X^+$ on the ruled
surface $X=P(L \oplus L_0)$ such that the real part of $X$ consists
of $t+k$ tori. Choose $k$ of these tori and make an elementary
transformation on each of them, that is the composition of the blowing
up at one point and the blowing down of the strict transform of the
fiber passing through this point. The result is still a real ruled surface of
base $(B,c_B)$ and the real part of this ruled surface consists of $t$
tori and $k$ Klein bottles, hence the result. $\square$

\subsection{The deformation's theorem}
\label{subsectiondeformation}

Let $\Delta \subset \C$ be the Poincar\'e's disk equipped with the complex
conjugation $conj$. A {\it real deformation} of surfaces is a proper
holomorphic submersion $\pi : Y \to \Delta$ where $(Y,c_Y)$ is a real
analytic manifold of dimension $3$ and $\pi$ satisfies $\pi \circ c_Y =
conj \circ \pi$. When $t \in ]-1 , 1 [ \in \Delta$, the fibers $Y_t =
\pi^{-1} (t)$ are invariant under $c_Y$ and are then compact real
analytic surfaces. Two real analytic surfaces $X'$ and $X''$ are said
to be 
{\it in the same deformation class} if there exists a chain $X'=X_0,
\dots , X_k=X''$ of compact real
analytic surfaces such that for every $i \in \{ 0 , \dots , k-1 \}$, the
surfaces $X_i$ and $X_{i+1}$ are isomorphic to some real fibers of a
real deformation.

\begin{prop}
The topological type of a real non-rational ruled surface is invariant under 
deformation.
\end{prop}

{\bf Proof :}

Let $(X,c_X) \to (B,c_B)$ be a real ruled surface of topological type
$(t,k,g,\mu , \epsilon)$ with $g \geq 1$. Let $\pi : Y \to \Delta$ be a real
deformation of surfaces such that $(Y_0 , c_Y|_{Y_0}) = (X,c_X)$. Then
every fiber of $\pi$ is a ruled surface with base of genus
$g$ (see \cite{BPV} for instance). 
Now since the deformation is trivial from the differentiable
point of view, the topology of the real part and the topology of the
involution on the base are invariant under deformation, hence the result. 
$\square$\\

For the sake of completeness, let us recall the following well known result,
see \cite{KhDg2} or \cite{Kh} :

\begin{theo}
\label{theodefrat}
There are four deformation classes of real structures on rational ruled
surfaces, one for which the real part is a torus, one for which the real 
part is a sphere and two for which the real part is empty. These two
later have non-homeomorphic quotients. $\square$
\end{theo}

Remember that the real structure for which the real 
part is a sphere is very special. It only exists on $\C P^1 \times \C P^1$
and is fibered over no real structure on the base $\C P^1$. This comes from
the existence of two rulings on $\C P^1 \times \C P^1$ and the involution 
$(x,y) \mapsto 
(y,x)$ reversing them. This is the main reason why we do not include
the case of rational ruled surfaces in theorem \ref{theoremdeformation}.

\begin{theo}
\label{theoremdeformation}
Two real non-rational ruled surfaces are in the same
deformation class if and only if they have the same topological type 
$(t,k,g, \mu
,\epsilon)$, except when $\mu=0$. There are two deformation classes of
real non-rational ruled surfaces of topological type $(0,0,g,0,0)$. For one 
such
class of ruled surfaces $(X, c_X)$, the quotient $X' = X/c_X$ is spin,
for the other one it is not.
\end{theo}

Using the terminology introduced in \cite{KhDg2}, this means that real
ruled surfaces are quasi-simple. The definition of the topological
type of a real ruled surface is given in \S
\ref{subsectiontoptype}. Note that every allowable quintuple is the topological
type of a real ruled surface (see proposition \ref{proposexistence}). 

\begin{rem}
If $X = P(E)$ is a real non-rational ruled surface of topological type 
$(t,k,g, \mu
,\epsilon)$ with $t+k < \mu$ and $k \neq 0$, then $X$ is
not decomposable, whereas any other topological type is
realized by a decomposable real ruled surface. Remember also that the 
deformation classes of complex ruled surfaces are described by the genus
 of the base and by whether the surface is spin or not. Then, real structures
for which $k$ is even only exist on spin ruled surfaces and real structures
for which $k$ is odd only exist on non-spin ruled surfaces.
\end{rem}

Let us sketch the proof of theorem \ref{theoremdeformation} :

Let $(X , c_X)$ be a real ruled and non-decomposable
surface with base $(B, c_B)$. If $X$ admits a real holomorphic
section, then we will prove that $(X , c_X)$ is in the same
deformation class that a real decomposable ruled surface (see
proposition \ref{proprealholsect}). If $X$ does not admit a real holomorphic
section, then we will prove that there exists a complex line bundle $L
\in \Pic (B)$ satisfying $c_B^* (L) = L^*$, such that $(X , c_X)$ is 
in the same
deformation class that the surface obtained from $(P(L \oplus L_0),
c_X^\pm)$ after at most one elementary transformation on each
component of its real part (see proposition \ref{propelemtrans}).

After these two steps, it is possible to reduce the study of
deformation classes of real strucures on ruled surfaces to the study
of deformation classes of real strucures on decomposable ruled
surfaces. It suffices then to check the  theorem
\ref{theoremdeformation} for decomposable real ruled surfaces.

\begin{prop}
\label{proprealholsect}
Let $(X , c_X)$ be a real ruled surface of base $(B, c_B)$ which
admits a real holomorphic section. Then there exists a real
deformation $\pi : Y \to \Delta$ such that for every $t \in \R^* \cap
\Delta$, $(Y_t , c_Y|_{Y_t})$ is isomorphic to $(X , c_X)$ and such
that $(Y_0 , c_Y|_{Y_0})$ is isomorphic to $(P(L \oplus L_0),
c_L \oplus c_{L_0})$ where $L \in \Pic (B)$ and $c_L$ is a real
structure on $L$ which lifts $c_B$.
\end{prop}

The definition of a real deformation has been given in the begining of
\S \ref{subsectiondeformation}.\\

{\bf Proof :}

Let $E$ be a rank two complex vector bundle over $B$ such that
$X=P(E)$. The real holomorphic section of $X$ is given by a complex
sub-line bundle $M$ of $E$. Denote by $N$ the quotient line bundle
$E/M$ so that the bundle $E$ is an extension of $N$ by $M$. Let $\mu \in H^1
(B , M \otimes N^*)$ be the extension class of this bundle and let
$\mu^1$ be a $1$-cocycle with coefficients in the sheaf ${\cal O}_B ( M
\otimes N^*)$, defined on a covering ${\cal U}=(U_i)_{i \in I}$ of
$B$, realising the cohomology class $\mu \in H^1
(B , M \otimes N^*)$. The bundle $E$ is then obtained as the gluing of the
bundles $(M \oplus N)|_{U_i}$ by the gluing maps :

\begin{eqnarray*}
(M \oplus N)|_{U_i \cap U_j} & \to & (M \oplus N)|_{U_j \cap U_i} \\
(m,n) & \mapsto &  \left[ \begin{array}{cc}
1&\mu_{ij}\\
0&1
\end{array} \right]
\left( \begin{array}{c}
m\\
n
\end{array} \right)
=(m+\mu_{ij} n , n).
\end{eqnarray*}
We can assume that for every open set $U_i$ of ${\cal U}$, there
exists $\overline{\i} \in I$ such that $U_{\overline{\i}} = c_B (U_i)$
(add these open sets to ${\cal U}$ if not). We can also assume that
there exists $J \subset I$ such that the open sets $(U_i)_{i \in J}$
cover $B$ and such that the real structure $c_X : X|_{U_i} \to
X|_{U_{\overline{\i}}}$ lifts to an antiholomorphic map $E|_{U_i} \to
E|_{U_{\overline{\i}}}$ (take a refinement of ${\cal U}$ if
not). Since by hypothesis the section of $X$ associated to $M$ is
real, these antiholomorphic maps are of the form :
\begin{eqnarray*}
(M \oplus N)|_{U_i} & \to & (M \oplus N)|_{U_{\overline{\i}}} \\
(x,(m,n)) & \mapsto &  (c_B (x) , \left[ \begin{array}{cc}
a_i&b_i\\
0&d_i
\end{array} \right]
\left( \begin{array}{c}
m\\
n
\end{array} \right) ),
\end{eqnarray*}
where $a_i$ (resp. $b_i$, resp. $d_i$) is an antiholomorphic morphism
$M|_{U_i} \to M|_{U_{\overline{\i}}}$ (resp. $N|_{U_i} \to
M|_{U_{\overline{\i}}}$, resp. $N|_{U_i} \to N|_{U_{\overline{\i}}}$)
which lifts $c_B$. Since $c_X$ is an involution, we have for every $i
\in J$, $a_{\overline{\i}} \circ a_i = d_{\overline{\i}} \circ d_i
\in {\cal O}_B^*|_{U_i}$ and $a_{\overline{\i}} \circ b_i + b_{\overline{\i}} \circ d_i
=0 \in {\cal O}_B (N^* \otimes M)|_{U_i}$. Moreover, for $i,j \in J$
such that $U_i \cap U_j \neq \emptyset$, the gluing conditions are the
following : $a_i = \lambda a_j$, $d_i = \lambda d_j$ and $b_i +
\mu_{\overline{\i} \overline{\j}} \circ d_i = \lambda (a_j \circ
\mu_{ij} + b_j)$ where $\lambda \in {\cal O}_B^*|_{U_i \cap U_j}$.

Now let $Y$ be the complex analytic manifold of dimension three
defined as the gluing of the charts $\C \times P(M \oplus N)|_{U_i}$, $i \in
J$, with change of charts given by the maps :
\begin{eqnarray*}
\C \times P(M \oplus N)|_{U_i} & \to & \C \times P(M \oplus N)|_{U_j} \\
(t,x,(m:n)) & \mapsto &  (t,x,\left[ \begin{array}{cc}
1&t \mu_{ij}\\
0&1
\end{array} \right]
\left( \begin{array}{c}
m\\
n
\end{array} \right))
=(t,x, (m+t \mu_{ij} n : n)).
\end{eqnarray*}
The projection on the first coordinate defines a holomorphic
submersion $\pi : Y \to \C$. The surface $\pi^{-1} (0)$ is isomorphic
to the decomposable ruled surface $P(M \oplus N)$, whereas, as soon as
$t \in \C^*$, the fiber $Y_t = \pi^{-1} (t)$ is isomorphic to the
ruled surface $X=P(E)$. Such an isomorphism $\psi_t : Y_t \to X$ is
given in the charts $P(M \oplus N)|_{U_i}$, $i \in J$, by :
\begin{eqnarray*}
P(M \oplus N)|_{U_i} & \to & P(M \oplus N)|_{U_j} \\
(x,(m:n)) & \mapsto &  (x, (m : t n)).
\end{eqnarray*}
Denote by $c_Y$ the real structure on $Y$ defined on charts $\C \times
P(M \oplus N)|_{U_i}$ by :
\begin{eqnarray*}
\C \times P(M \oplus N)|_{U_i} & \to & \C \times P(M \oplus N)|_{U_{\overline{\i}}} \\
(t,x,(m:n)) & \mapsto &  (\overline{t},c_B (x),\left[ \begin{array}{cc}
a_i&\overline{t} b_i\\
0&d_i
\end{array} \right]
\left( \begin{array}{c}
m\\
n
\end{array} \right)).
\end{eqnarray*}
This real structure satisfies $\pi \circ c_Y = conj \circ \pi$ where
$conj$ is the complex conjugation on $\C$. Moreover, when $t \in
\R^*$, $\phi_t$ gives an isomorphism between the real ruled surfaces
$(Y_t , c_Y|_{Y_t})$ and $(X , c_X)$. Hence, the restriction of $\pi :
Y \to \C$ over $\Delta \subset \C$ is a real deformation which
satisfies proposition \ref{proprealholsect}. $\square$

\begin{prop}
\label{propelemtrans}
Let $(X , c_X)$ be a real ruled surface of base $(B, c_B)$, which
does not admit any real holomorphic section. Then, there exists $L \in
\Pic (B)$ satisfying $c_B^* (L)=L^*$ and a ruled surface $(X' ,
c_{X'})$ obtained from $(P(L \oplus L_0) , c_X^\pm)$ after at most one
elementary transformation on each of its real components, such that
$(X , c_X)$ and $(X' , c_{X'})$ are in the same deformation class.
\end{prop}
Remember that an {\it elementary transformation} on the ruled surface
$X$ is by definition the composition of a blowing up of $X$ at one
point and the blowing down of the strict transform of the fiber
passing through this point.

\begin{lemma}
\label{lemmaelemtransfsect}
Let $X=P(L \oplus L_0)$ be a decomposable ruled surface of base
$B$. Let $s : B \to X$ be the section defined by $L$ and $D$ be a
divisor associated to $L$. Then the ruled surface obtained from $X$
after an elementary transformation at the point $s (x)$, $x \in B$, is
the surface $P(L(x) \oplus L_0)$ where $L(x)$ is the complex line
bundle associated to the divisor $D+x$. $\square$
\end{lemma}

\begin{lemma}
\label{lemmaveryample}
Let $(X , c_X)$ be a real ruled surface of base $(B, c_B)$, which
does not admit any real holomorphic section. Then $X$ has a very ample
holomorphic section $S$ which is transversal to its image under $c_X$.
\end{lemma}

{\bf Proof :}

Let us first construct a very ample section on $X$. Let $E$ be a rank
two complex vector bundle over $B$ such that $X=P(E)$, and let $A$ be
an ample line bundle over $B$. Then by definition, for sufficiently large $n$,
the bundle $E^* \otimes A^n$ is generated by its global sections. Choosing
$N$ such global sections, it provides a surjective morphism of bundles
$B \times \C^N \to E^* \otimes A^n$. This induces an injective morphism
between the dual bundles $E \otimes (A^*)^n \to B \times \C^N$ and thus
an embedding $X \to B \times \C P^{N-1}$. Fixing an embedding 
$B \to \C P^3$, we deduce an embedding $X \to \C P^3 \times \C P^{N-1}$.
Finally, combining this with Segre embedding, we obtain an embedding
$X \to \C P^{4N-1}$ associated to a very ample linear system of sections 
on $X$.

Now, let us prove that in this linear system, there exists a smooth section $S$
transversal to $c_X (S)$. From Bertini's theorem (see \cite{Hart}, theorem $8.18$)
there exists, in this linear system, a smooth
section $S$ associated to a hyperplane $H$ of
$\C P^{4N-1}$ transversal to $X$. By hypothesis, $S$ cannot be real, so that
the intersection $c_X (S) \cap S$ consists of a finite number of points. 
We will prove that after a small perturbation of $H$, this intersection
can be assumed transversal. Indeed, let $x \in c_X (S) \cap S$. If $x \in \R X$,
the intersection of $H$ with $T_x X$ is a line, which is the tangent of $S$ at $X$.
The section $S$ is transverse to $c_X (S)$ at $x$ if and only if this line is
not fixed by the differential $d_x c_X$. Since the fixed point set of this involution 
is of half 
dimension, the intersection of $S$ and $c_X (S)$ at $x$
can be made transversal after a small perturbation of $H$, keeping the
intersection point $x$. Now, if $x \notin \R X$, then since the section $S$ is
smooth, the points $x$ and $c_X (x)$ belong to two different fibers of $X$ and in
particular to non-real ones.
Suppose that the line $D_x \subset \C P^{4N-1}$ joining them is transversal
to both the planes $T_x X$ and $T_{c_X (x)} X$. Then there exists a pencil 
of hyperplanes of $\C P^{4N-1}$ containing $H$ and parametrised both by the lines of 
$T_x X \subset \C P^{4N-1}$ and the lines of $T_{c_X (x)} X \subset \C P^{4N-1}$.
This means that each line of $T_{x} X$ passing through $x$, and similarly
each line of $T_{c_X (x)}$ passing through $c_X (x)$, is contained in one and
only one hyperplane of this pencil. Also, this pencil contains no other hyperplane.
$$\vcenter
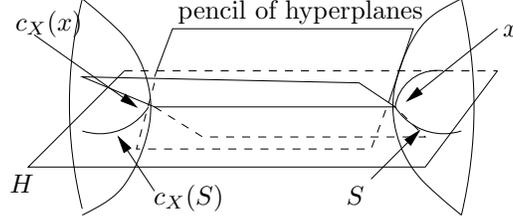
{\hbox{\begin{picture}(0,0)%
\epsfig{file=realstr2.pstex}%
\end{picture}%
\setlength{\unitlength}{3315sp}%
\begingroup\makeatletter\ifx\SetFigFont\undefined%
\gdef\SetFigFont#1#2#3#4#5{%
  \reset@font\fontsize{#1}{#2pt}%
  \fontfamily{#3}\fontseries{#4}\fontshape{#5}%
  \selectfont}%
\fi\endgroup%
\begin{picture}(3690,1644)(1711,-4168)
\put(5401,-2806){\makebox(0,0)[lb]{\smash{\SetFigFont{10}{12.0}{\rmdefault}{\mddefault}{\updefault}$x$}}}
\put(1756,-2761){\makebox(0,0)[lb]{\smash{\SetFigFont{10}{12.0}{\rmdefault}{\mddefault}{\updefault}$c_X (x)$}}}
\put(2791,-4066){\makebox(0,0)[lb]{\smash{\SetFigFont{10}{12.0}{\rmdefault}{\mddefault}{\updefault}$c_X (S)$}}}
\put(1711,-3976){\makebox(0,0)[lb]{\smash{\SetFigFont{10}{12.0}{\rmdefault}{\mddefault}{\updefault}$H$}}}
\put(4231,-4066){\makebox(0,0)[lb]{\smash{\SetFigFont{10}{12.0}{\rmdefault}{\mddefault}{\updefault}$S$}}}
\put(2971,-2671){\makebox(0,0)[lb]{\smash{\SetFigFont{10}{12.0}{\rmdefault}{\mddefault}{\updefault}pencil of hyperplanes}}}
\end{picture}
}}$$
This pencil thus provides us with a holomorphic identification between the projective
lines $P(T_{x} X)$ and $P(T_{c_X (x)} X)$. Under this identification, the 
differential $d_x c_X$ reads as an anti-holomorphic involution of $T_{x} X$ and once
 more,
the section $S$ is transversal to $c_X (S)$ at $x$ if and only its tangent line is
not fixed by this involution $d_x c_X$. This can always be garanted after a small
perturbation of $H$. Since small perturbations do not perturb the transversality of
transversal points, this process strictly increases
the number of transversal points between $S$ and $c_X (S)$ and so gives the result after
a finite number of steps. It thus only remains to prove that the line $D_x$ can 
indeed be assumed transverse to both the planes $T_x X$ and $T_{c_X (x)} X$, after a
small perturbation of $H$ if necessary.

For this, note that the embedding $B \to \C P^3$ can be chosen real. The set of points
of $B$ whose tangent is not a real line of $\C P^3$ is a then dense open subset 
$U \subset B$
(for the usual topology, not the Zariski's one), invariant under $c_B$. The set $U$
 is in fact
the complementary of the real part of the dual curve. Let $x \in X$ be a point such that
$y=p(x) \in U$ where $p$ is the projection $X \to B$. Since the line joining $y$ to 
$c_B (y)$ is real, it is not tangent to $B$ at $y$ and $c_B (y)$. Let $H_1$ be a 
hyperplane
of $\C P^3$ passing through $y$ and $c_B (y)$ and transverse to $B$. Then $H_1 \times 
\C P^{N-1}$ is transverse to $X$ in $\C P^3 \times 
\C P^{N-1}$. Let $H_2$ be a hyperplane of $\C P^{N-1}$ such that $\C P^3 \times H_2$ 
does contain neither $x$ nor $c_X (x)$. Then the divisor $(H_1 \times 
\C P^{N-1}) + (\C P^3 \times H_2)$ is associated to a hyperplane $H_0$ of $\C P^{4N-1}$,
which contains both $x$ and $c_X (x)$ and which is transverse to $X$ at these points.
Then $H_0$ contains the line $D_x$ and since by construction it also contains the fibers 
through $x$ and $c_X (x)$, its transversality with $X$ at $x$ and $c_X (x)$ implies
the one of $D_x$. Hence for any point $x$ belonging to the open set $p^{-1} (U)$
of $X$, the line $D_x$ is transverse to $X$ at $x$ and $c_X (x)$. Since it is not hard
to observe that any non-real intersection point of $S$ and $c_X (S)$ can be moved
to $p^{-1} (U)$ after a small perturbation of $H$, this completes to proof of
lemma \ref{lemmaveryample}. $\square$ \\

{\bf Proof of proposition \ref{propelemtrans} :}

Let $S \subset X$ be a very ample holomorphic smooth section, transverse
to its image under $c_X$. Such a section is given by lemma \ref{lemmaveryample}.
The set $c_X (S) \cap S$ is finite and invariant
under $c_X$. Denote by $X_1$ the 
ruled surface obtained from $X$ after
an elementary transformation on every point of this set. Since
it is invariant under $c_X$, the real structure $c_X$
induces a real structure $c_{X_1}$ on $X_1$. Moreover, the strict
transform $S_1$ of $S$ satisfies $c_{X_1} (S_1) \cap S_1 =
\emptyset$. Thus $X_1$ is a decomposable ruled surface, and $c_{X_1}$
exchanges the two holomorphic sections $S_1$ and $c_X (S_1)$. The
inverse of an elementary transformation is still an elementary
transformation, so we deduce that $(X, c_X)$ is obtained from the real
decomposable ruled surface $(X_1 , c_{X_1})$ after performing
elementary transformations on points $\{ x_1 , \dots , x_k , y_1 ,
\dots , y_l , \overline{y}_1 ,
\dots , \overline{y}_l \}$ where $c_{X_1} (x_i) = x_i$ and $c_{X_1} (y_j)
=\overline{y}_j$. Note that all the points $\{ x_1 , \dots , x_k , y_1 ,
\dots , y_l , \overline{y}_1 ,
\dots , \overline{y}_l \}$ belong to different fibers of $X_1$. It remains 
to see that this number of points can be
reduced to one at most for each component of $\R X_1$, changing the
decomposable real ruled surface $X_1$ if necessary.

For every $j \in \{1 , \dots , l\}$, choose a piecewise analytic path
$y_j (t)$, $t \in [0,1]$, such that $y_j (0) = y_j$, $y_j (1) \in S_1$
and  $p(y_j (t))$ is constant, which means that $y_j (t)$ stays in a
same fiber of $X_1$. Let
$\overline{y}_j (t) = c_{X_1} (y_j (t))$ and denote by $X_2$ the ruled
surface obtained from $X_1$ after elementary transformations in the
points $y_1 (1) , \dots , y_l (1) , \overline{y}_1 (1) , \dots ,
\overline{y}_l (1)$. The real structure  $c_{X_1}$ induces a real
structure  $c_{X_2}$ on $X_2$. The surface $(X_2 , c_{X_2})$ is in the
same deformation class that $(X_1 , c_{X_1})$. Moreover, $X_2$ is also
a decomposable ruled surface. Indeed, the strict transform $S_2$ of
$S_1$ is a holomorphic section of $X_2$ satisfying $c_{X_2} (S_2) \cap
S_2 = \emptyset$. Thus $(X, c_X)$ is in the same deformation class
that the surface obtained from the real
decomposable ruled surface $(X_2 , c_{X_2})$ after performing
elementary transformations on the strict transforms of the points 
$x_1 , \dots , x_k \in \R X_1$, still denoted by 
$x_1 , \dots , x_k \in \R X_2$. Now for
each  pair of points $x_1 , x_2$ lying in a
same connected component of $\R X_2$, we can make the elementary transformation
on the point $x_2$. Then, the image of the fiber passing through $x_2$ is
a real point $x'_2$ in the new surface $X'_2$ obtained. So we can choose an
analytic path from $x_1$ to $x'_2$ in the real part of $X'_2$ and we
deduce that the surface obtained from $X_2$ after making the elementary 
transformations
on the points $x_1 , x_2$ is in the same deformation class that the one
 obtained from $X'_2$ after an elementary transformation on $x'_2$, which is 
$X_2$ itself. Hence each pair of points lying in a
same connected component of $\R X_2$ can be removed and so $(X, c_X)$ is in 
the same deformation class
that the surface obtained from the real
decomposable ruled surface $(X_2 , c_{X_2})$ after performing at most one
elementary transformation on each of its real components.
Since $c_{X_2}$ exchanges two disjoint holomorphic sections of
$X_2$, it follows from theorem \ref{theoremrealstruct} that $(X_2 ,c_{X_2})$ 
is of the form $(P(L \oplus L_0) , c_X^+)$ where $L \in \Pic
(B)$ and $c_B^* (L) = L^*$. $\square$

\begin{lemma}
\label{lemmespin}
Let $g \geq 1$ be an odd integer and $(B,c_B)$ be a smooth compact
irreducible real algebraic curve of genus $g$ and empty real
part. Let $L$ be a complex line bundle over $B$ satisfying
$c_B^* (L) = L^*$. Then the real ruled surfaces $(P(L \oplus L_0), c_{X}^+)$
 and $(P(L \oplus L_0), c_{X}^-)$ are in the same deformation class.
\end{lemma}
(In lemma \ref{lemmespin}, the real structures $c_{X}^+$ and $c_{X}^-$ on
$X= P(L \oplus L_0)$ are those given by proposition \ref{propRL*}.)\\

{\bf Proof :}

Without changing the deformation class of $X= P(L \oplus L_0)$,
we can assume that the base of this surface is the real algebraic curve
$(B,c_B)$ given by corollary \ref{corspin}. Then, if $L$ belong to the same 
real component of $(\Jac (B) , -c_B^*)$ that the bundle given by corollary 
\ref{corspin}, we can assume, without changing the deformation class of 
$X= P(L \oplus L_0)$, that $L$ is exactly this bundle.
In that case, the result comes from corollary \ref{corspin}.

Let $X= P(L \oplus L_0)$ be the ruled surface given by corollary \ref{corspin},
and $\Phi : X \to X$ be the automorphism conjugating $c_{X}^+$ and $c_{X}^-$.
Let $x_1$ be a point on the section of $X$ associated to $L$ and
$y_1 = c_{X}^+ (x_1) = c_{X}^- (x_1)$. Let $x_2 = \Phi (x_1)$ and
$y_2 = \Phi (y_1) = c_{X}^+ (x_2) = c_{X}^- (x_2)$. Denote by $Y_1$ 
(resp. $Y_2$) the ruled surface obtained from $X$ after one elementary
transformation on the points $x_1$ and $y_1$ (resp. $x_2$ and $y_2$). Then
the real structures $c_{X}^+$ and $c_{X}^-$ lift to the real structures 
$c_{Y_1}^\pm$ (resp. $c_{Y_2}^\pm$) on $Y_1$ (resp. $Y_2$), and 
$\Phi$ lifts to a biholomorphism $\Psi : Y_1 \to Y_2$ such that
$c_{Y_1}^+ = \Psi^{-1} \circ c_{Y_2}^- \circ \Psi$ and 
$c_{Y_1}^- = \Psi^{-1} \circ c_{Y_2}^+ \circ \Psi$. But the real ruled
surface $(Y_1 , c_{Y_1}^-)$ is in the same deformation class that
$(Y_2 , c_{Y_2}^-)$. Indeed, it suffices to choose an analytic path $x_t$
linking $x_1$ to $x_2$ in the section of $X$ associated to $L$ and to 
consider the surfaces $(Y_t ,c_{Y_t}^-)$ obtained from $(X, c_{X}^-)$
after an elementary transformation on the points $x_t$ and $c_{X}^- (x_t)$.

Hence the real ruled surfaces $(Y_1 , c_{Y_1}^-)$ and $(Y_1 , c_{Y_1}^+)$
are in the same deformation class. To conclude, it remains to see that
they do not come from the same connected component of $(\Jac (B) , -c_B^*)$ 
that $(X, c_{X}^\pm)$. This follows from the fact that the quotients
$Y_1 / c_{Y_1}^\pm$ and $X / c_{X}^\pm$ are not homeomorphic. Indeed,
these two quotients are sphere bundles over the non-orientable surface 
$B' = B / c_B$. But $Y_1 / c_{Y_1}^\pm$ is obtained from $X / c_{X}^\pm$
after one elementary transformation in one point. Thus one of these
two quotient is spin, and one is not. Hence the result. $\square$ \\

{\bf Proof of theorem \ref{theoremdeformation} :}

Let $(X_1 ,c_{X_1})$ and $(X_2 ,c_{X_2})$ be two real non-rational ruled 
surfaces
of bases $(B_1 ,c_{B_1})$ and $(B_2 ,c_{B_2})$ respectively, which have the 
same
topological type $(t,k,g,\mu , \epsilon)$. We have to prove that they
are in the same deformation class, as soon as $\mu \neq 0$.

Let us first consider the case of decomposable ruled surfaces, that is
let us assume that $X_1$ and $X_2$ are decomposable. If $t+k < \mu$,
it follows from theorem \ref{theoremrealstruct} that $X_1 = P(L_1
\oplus L_0)$ (resp. $X_2 = P(L_2
\oplus L_0)$), where $L_1 \in \Pic (B_1)$ (resp. $L_2 \in \Pic (B_2)$)
and $c_{B_1}^* (L_1) = L_1^*$ (resp. $c_{B_2}^* (L_1) = L_2^*$). Moreover,
it follows from proposition \ref{propRL*} that in this case $k=0$. The
partition ${\cal
P}_1$ (resp. ${\cal
P}_2$) in two elements of $\R B_1$ (resp. $\R B_2$) associated to
$L_1$ (resp. $L_2$) consists of one element containing $t$ components
of $\R B_1$ (resp. $\R B_2$) and one element containing $\mu - t$ components
of $\R B_1$ (resp. $\R B_2$) (see \S \ref{subsectionpartition} for the
definition of the partition). Since $(B_1 ,c_{B_1})$ and $(B_2
,c_{B_2})$ have same
topological type $(g,\mu , \epsilon)$, there exists a piecewise
analytic path of smooth real algebraic curves connecting them (see
\cite{Nat2}). Moreover, this path can be chosen such that the $t$
components of $\R B_2$, which form an element of the partition ${\cal
P}_2$, deform into the $t$
components of $\R B_1$ which form an element of the partition ${\cal
P}_1$. This follows from the presentation in \cite{Nat2}
of a real algebraic curve as the gluing of a Riemann surface with
boundary with its conjugate, the gluing maps being either identity or
antipodal. Thus $(X_2 ,c_{X_2})$ is in the same deformation class that
a ruled surface $(\widetilde{X}_2 ,c_{\widetilde{X}_2})$ of base $(B_1
,c_{B_1})$. Moreover, 
$\widetilde{X}_2 = P(\widetilde{L}_2 \oplus L_0)$ where
$\widetilde{L}_2 \in \Pic (B_1)$, $c_{B_1}^* (\widetilde{L}_2) =
\widetilde{L}_2^*$ and the partitions associated to $\widetilde{L}_2$
and $L_1$ are the same. From lemma \ref{lemmapartition} follows that 
$\widetilde{L}_2$
and $L_1$ are in the same component of the real part of $(\Jac (B_1) ,
-c_{B_1}^*)$ and hence the surfaces $(\widetilde{X}_2 ,
c_{\widetilde{X}_2})$ and $(X_1 ,c_{X_1})$ are in the same deformation
class.

If $t+k = \mu$, it follows from theorem \ref{theoremrealstruct} that 
$X_1 = P(L_1
\oplus L_0)$ (resp. $X_2 = P(L_2
\oplus L_0)$), where $L_1 \in \Pic (B_1)$ (resp. $L_2 \in \Pic (B_2)$)
and either $c_{B_1}^* (L_1) = L_1^*$ (resp. $c_{B_2}^* (L_1) =
L_2^*$), or $c_{B_1}^* (L_1) = L_1$ (resp. $c_{B_2}^* (L_1) =
L_2$). In the first case, $L_1$ (resp $L_2$) is in the same component
of the real part of $(\Jac (B_1) , -c_{B_1}^*)$ (resp. $(\Jac (B_2) ,
-c_{B_2}^*)$) that $L_0$, since $t+k = \mu$. Thus $(X_1 ,c_{X_1})$
(resp. $(X_2 ,c_{X_2})$) is in the same deformation class that $(B_1
\times \C P^1 , c_X^\pm)$ (resp. $(B_2
\times \C P^1 , c_X^\pm)$). Moreover, when $\mu \neq 0$, only one of
the two real structures $c_X^\pm$, say $c_X^+$, satisfies $t+k =
\mu$. In the second case, denote by $D_+ - D_-$ a divisor associated
to $L_1$, where $D_+$, $D_-$ are positive divisors and invariant under
$c_{B_1}$. Then $X_1 = P(L_{D_+} \oplus L_{D_-})$ and $c_{X_1} =
c_{L_{D_+}} \oplus c_{L_{D_-}}$. Thus, it follows from lemma
\ref{lemmaelemtransfsect} that $(X_1 ,c_{X_1})$ is obtained from $(B_1
\times \C P^1 , c_{L_0} \oplus c_{L_0})$ after performing elementary
transformations on the points of the section associated to $L_{D_+}$
(resp. $L_{D_-}$) over the locus of $D_+ \in B_1$ (resp. $D_- \in
B_1$). Without changing the deformation
class of the surface, we can assume that the elementary
transformations are only done on real points of $(B_1
\times \C P^1 , c_{L_0} \oplus c_{L_0})$ with at most one on each of
its real components. Indeed, the extra real points can be removed as in 
proposition \ref{propelemtrans} and every couple of conjugated imaginary points
can be moved to real points following a standard deformation : embedd
the disk $(\Delta , conj)$ in a real section of $X$, and for every 
$t \in \Delta$, denote by $Y_t$ the surface obtained from $X$ after an
elementary transformation on the points $t$ and $-t$ in $\Delta$ (we
still denote by $\Delta$ its image in $X$ by the chosen embedding). The 
dimension $3$ complex manifold $Y$ obtained gets two real structures,
one which lifts $conj$ in $\Delta$ and one which lifts $-conj$. This
thus define two real deformations of ruled surfaces and shows that the
real ruled surfaces obtained from $X$ after making elementary transformations
on the points $\pm \frac{1}{2} \in \Delta$ or $\pm \frac{i}{2} \in \Delta$
are in the same deformation class. Hence, without changing the deformation
class of the surface $(X_1 ,c_{X_1})$, we can assume that the elementary
transformations are done only on real points of $(B_1
\times \C P^1 , c_{L_0} \oplus c_{L_0})$ with at most one on each of
its real components. The total number of such elementary
transformations is then $k$ since the topological type of $(X_1
,c_{X_1})$ is $(t,k,g,\mu , \epsilon)$. If $X_1$ and $X_2$ are two
such surfaces,  there exists a piecewise
analytic path of smooth real algebraic curves connecting $(B_1 ,c_{B_1})$ and 
$(B_2
,c_{B_2})$, such that the $k$ components of $\R B_2$ over which are
done the elementary transformations deform on the $k$ components of
$\R B_1$ over which are 
done the elementary transformations. Hence in both cases, $(X_1 ,c_{X_1})$ 
and $(X_2
,c_{X_2})$ are in the same deformation class. Since the real
structures $c_X^+$ and $c_{L_0} \oplus c_{L_0}$ are conjugated on $B_1
\times \C P^1$, which follows from theorem \ref{theoremrealstruct} for
instance, we deduce that the real decomposable ruled surfaces $(X_1
,c_{X_1})$ and $(X_2 ,c_{X_2})$ are in the same deformation class if
and only if they have the same topological type $(t,k,g,\mu ,
\epsilon)$, except when $\mu = 0$. In that case, if $g$ is even, it
follows from proposition \ref{propjac} that the same method as before
leads to the fact that $(X_1
,c_{X_1})$ and $(X_2 ,c_{X_2})$ are in the same deformation class that
$(B \times \C P^1 , c_X^+)$ or $(B \times \C P^1 , c_X^-)$. But the
quotient $(B \times \C P^1) / c_X^+$ is spin and $(B \times \C P^1) /
c_X^-$ is not, so the surfaces $(B \times \C P^1 , c_X^+)$ and $(B
\times \C P^1 , c_X^-)$ are not in the same deformation class. If $g$ is odd, 
it
follows from proposition \ref{propjac} that the same method as before
leads to the fact that $(X_1
,c_{X_1})$ and $(X_2 ,c_{X_2})$ are in the same deformation class that
$(P(L \oplus L_0) , c_X^\pm)$, where $L$ belongs to one of the two
components of the real part of $( Jac (B) , -c_B^*)$. But it follows
from lemma \ref{lemmespin} that $(P(L \oplus L_0) , c_X^+)$ and
$(P(L \oplus L_0) , c_X^-)$ are in a same deformation class. The
result follows from the fact that $(B \times \C P^1) / c_X^\pm$ is
spin and  $P(L \oplus L_0) / c_X^\pm$ is not when $L$ is not in the
same component of the real part of $( Jac (B) , -c_B^*)$ that $L_0$.

Now let us prove the theorem in the general case, which means that we no
more assume that $X_1$ and $X_2$ are decomposable. From propositions
\ref{proprealholsect} and \ref{propelemtrans} follow that these
surfaces are either in the same deformation class that some real
decomposable ruled surfaces, or in the same deformation class that
some ruled surface obtained from a decomposable one of the form $(P(L
\oplus L_0) , c_X^\pm)$ after at most one
elementary transformation on each of its real components. In this
second case, we can assume that $L$ does not belong to the same
component of the real part of $( Jac (B) , -c_B^*)$ that $L_0$
(otherwise the surface can be deformed to a decomposable ruled
surface). Since the topological types of these surfaces are different
from those realized by decomposable ruled surface, we can assume that
either $X_1$ and $X_2$ are both decomposable, or that they are both
from this second class. In the first case, the theorem follows from
what we have already done. Let us assume we are in the second
case. Then there exists $L_1 \in \Pic (B_1)$ (resp. $L_2 \in \Pic (B_2)$)
such that $c_{B_1}^* (L_1) = L_1^*$ (resp. $c_{B_2}^* (L_2) = L_2^*$)
and $(X_1 , c_{X_1})$ is obtained from $(P(L_1 \oplus L_0) , c_X^+)$
after making $k$ elementary transformations in $k$ disjoint real
components. The surfaces $(P(L_1 \oplus L_0) , c_X^+)$ and $(P(L_2
\oplus L_0) , c_X^+)$ have same topological type $(t+k , 0 , g , \mu ,
\epsilon)$, with $\mu >0$. Thus they are in the same deformation
class. Moreover, in the same way as before, this deformation can be
chosen so that the $k$ marked real components of $(P(L_2
\oplus L_0) , c_X^+)$ deforms to the $k$ marked real components of
$(P(L_1 \oplus L_0) , c_X^+)$. It follows that $(X_1
,c_{X_1})$ and $(X_2 ,c_{X_2})$ are in the same deformation class. $\square$

\addcontentsline{toc}{part}{\hspace*{\indentation}Bibliographie}

\nocite{*}  

\bibliography{realstruct}
\bibliographystyle{abbrv}

\noindent Ecole Normale Sup\'erieure de Lyon\\
Unit\'e de Math\'ematiques pures et appliqu\'ees\\
$46$, all\'ee d'Italie\\
$69364$, Lyon C\'edex $07$
(FRANCE)\\
e-mail : {\tt jwelschi@umpa.ens-lyon.fr}

\end{document}